\documentclass{amsart}[12pt]
\usepackage{amssymb,latexsym}

\numberwithin{equation}{section}

\newtheorem{theorem}{Theorem}[section]

\newtheorem{proposition}[theorem]{Proposition}
\newtheorem{corollary}[theorem]{Corollary}
\newtheorem{remark}[theorem]{Remark}
\newtheorem{lemma}[theorem]{Lemma}
\newtheorem{example}[theorem]{Example}
\newtheorem{definition}[theorem]{Definition}

\def\ZZ{\mathbb{Z}}
\def\CC{\mathbb{C}}

\def\l{\ell}
\def\wnot {w_\mathrm{o}}

\def\ii{\mathbf{i}}

\newcommand{\mat}[4]{\left(\!\!\begin{array}{cc}
#1 & #2 \\
#3 & #4 \\
\end{array}\!\!\right)}

\begin{document}

\addtocounter{section}{-1}

\title{Noncommutative Double Bruhat cells and their factorizations}

\date{March 21, 2004}

\author{Arkady  Berenstein}

\author{Vladimir Retakh}

\address{Department of Mathematics, University of Oregon}
\email{arkadiy@math.uoregon.edu}


\address{Department of Mathematics, Rutgers University}
\email{vretakh@math.rutgers.edu}

\date{June 30, 2004}

\thanks{Research supported in part by the NSF (DMS) grant \# 0102382 (A.B) and
by the NSA grant \# MDA 9040310047 (V.R)}


\maketitle

\makeatletter
\renewcommand{\@evenhead}{\tiny \thepage \hfill
A.~BERENSTEIN and V.~RETAKH \hfill}

\renewcommand{\@oddhead}{\tiny \hfill
NONCOMMUTATIVE DOUBLE BRUHAT CELLS AND THEIR FACTORIZATIONS
 \hfill \thepage}
\makeatother

\tableofcontents



\maketitle

\section{Introduction}
This paper is a first attempt to generalize results of A. Berenstein,
S. Fomin and  A. Zelevinsky on total positivity of matrices over
commutative rings to matrices over noncommutative rings.

The classical theory of total positivity studies matrices whose
minors all are nonnegative. Motivated by a surprising connection
discovered by G. Lusztig \cite {l1, l2}
between total positivity of matrices and canonical bases for
quantum groups, A. Berenstein, S. Fomin and A. Zelevinsky in a
series of papers \cite {bfz, bz0, bz, fz} systematically investigated the problem of
total positivity from a representation-theoretic point of view.

In particular, they showed that a natural framework for the study of totally
positive matrices is provided by the decomposition of a reductive
group $G$ into the disjoint union of double Bruhat cells
$G^{u,v}=BuB\cap B_-vB_-$ where $B$ and $B_-$ are two opposite
Borel subgroups in $G$, and $u$ and $v$ belong to the Weyl group
$W$ of $G$.

According to \cite {bfz, bz, fz} there exist families of
birational parametrizations of $G^{u,v}$, one for each reduced
expression of the element $(u,v)$ in the Coxeter group $W\times
W$. Every such parametrization can be thought of as a system of
local coordinates in $G^{u,v}$. Such coordinates are called {\it the
factorization parameters} associated to the reduced expression of
$(u,v)$. The coordinates are obtained by expressing a generic
element $x\in G$ as an element of the maximal torus $H=B\cap B_-$
multiplied by the product of elements of various one-parameter
subgroups in $G$ associated with simple roots and their
negatives; the reduced expression  prescribes the order of factors
in this product. An explicit formula for these factorization
parameters as rational functions on the double Bruhat cell
$G^{u,v}$ was given.

As we said, Berenstein, Fomin and Zelevinsky came to
factorization parameters (first, for $GL_n$ and then for other
classical groups)  from representation theory. For the
noncommutative case our program is to go into opposite direction:
from factorization parameters for $GL_n$ to ``total positivity'',
canonical bases  and representations. This paper is a beginning
of the program.

For $G=GL_n(F)$ where $F$ is a field of characteristic zero, the
explicit formulas for factorization parameters are given through
the classical determinant calculus. As a first step toward
noncommutative representation theory and noncommutative total
positivity, we generalize here the results  from \cite {fz} and
\cite {bz} to $G=GL_n({\mathcal F})$ where ${\mathcal F}$ is a (noncommutative)
skew field by using the {\it quasideterminantal} calculus of matrices
over (noncommutative) rings introduced by I. Gelfand and V. Retakh
\cite{gr1, gr2, gr3, ggrw}.

The noncommutative point of view has many advantages. Let $\wnot \in
W$ be the element of the maximal length. In the commutative case
the factorization parameters for $x\in G^{u,v}$, $G=GL_n$,
$u=id$, $v=\wnot $ are given as ratios $ab/cd$ or $a/b$ where $a,b,c,d$
are minors of matrix $x$ (see \cite {bfz}). In the noncommutative
case, {\it for any} $u$ and $v=\wnot $, the factorization parameters
can be written as $f^{-1}g$ where $f,g$ are quasiminors for
matrix $x$. The paper contains other noncommutative formulas and
constructions for $GL_n$ that are new even in the commutative case.

Our results confirm the Gelfand principle: noncommutative algebra
(properly understood) is simpler than its commutative counterpart.

The paper is organized as follows.

In Section \ref {sec:quasideterminants} we recall some facts
about quasideterminants and introduce our main tool - positive
quasiminors $\Delta^i_{u,v}$. In Section \ref{sec:basic
factorizations} we study basic factorizations in $GL_n$ and its
Borel subgroup. Section \ref{sec:examples} contains examples of
such factorizations. Section \ref{sec:double bruhat cells}
section is central for the paper. First, we introduce
``noncommutative $SL_2$-subgroups'' in $GL_n$. For a generic
matrix $x$ we define special quasiminors $\Delta ^i_{u,v}(x)$,
where $u,v\in W$ and show that they satisfy certain ``Pl\"ucker
relations". We note that $\Delta^i_{u,v}(x)$ is always positive
for positive real matrices $x$. Section \ref{sec:double bruhat
cells} also contains the main result: it gives formulas for
factorization coordinates for reduced double Bruhat cells. For a
matrix $x\in G^{u,v}$ these coordinates are written as products
of quasiminors $\Delta^i_{s,t}(y)$ where the matrix $y$ is the so called
noncommutative twist of $x$. In Section \ref{sec:other
factorizations} we study relations between quasiminors of $x\in
G^{u, \wnot }$ and the corresponding twisted matrix. In this case the
quasiminors $\Delta_{\cdot,\cdot}^i(y)$ in the main theorem can
be replaced by quasiminors $\Delta_{\cdot,\cdot}^i(x)$. Studying
twisted matrices is an important problem by itself and we present
several approaches to computations of such matrices. These
results are new even in the commutative case.

\section{Quasideterminants and Quasiminors}

\label{sec:quasideterminants}

A notion of quasideterminants for matrices over a noncommutative
ring was introduced in \cite{gr1, gr2} and developed in \cite{gr3}. It has been
effective in many areas (see, for the example, the survey article \cite {ggrw}).
Here we remind a few facts about quasideterminants which will be used in this
paper.

\subsection{Definition of quasideterminants}
Let $A=(a_{ij}), i\in I, j\in J$ be a matrix of order $n$ over a ring $R$.
Construct the following submatrices of A:
submatrix $A^{ij}$, $i\in I$, $j\in J$ obtained from $A$ by deleting its $i$-th
row and $j$-th column;
row submatrix $r_i$ obtained from $i$-th row of $A$ by deleting the
element $a_{ij}$;
column submatrix $c_j$ obtained from $j$-th column of $A$ by deleting the
element $a_{ij}$.

\begin{definition}
{\rm If $n=1$  the quasideterminant $|A|_{ij}$
equals to $a_{ij}$. If $n>1$   the quasideterminant $|A|_{ij}$ is
defined if the submatrix $A^{ij}$ is invertible over the ring
$R$. In this case one has
$$
|A|_{ij}=a_{ij}-r_i(A^{ij})^{-1}c_j.
$$
}
\end{definition}

For a generic matrix $A$ over a skew field ${\mathcal F}$, one has

$$|A|_{ij} = a_{ij} -\sum
a_{iq}|A^{ij}|^{-1}_{pq} a_{pj}.
$$
Here the sum is taken over all $p\in I\smallsetminus\{i\},q\in
J\smallsetminus\{j\}$.

If $A$ is an $n\times n$-matrix there exist up to $n^2$
quasideterminants of $A$.

By definition, an  $r$-quasiminor of a square matrix $A$ is a
quasideterminant of an $r\times r$-submatrix of $A$.

Sometimes it is convenient to adopt a more graphic notation for
the quasideterminant $|A|_{pq}$ by boxing the element $a_{pq}$.
For $A=(a_{ij})$, $i,j=1,\dots , n$, we write
$$
|A|_{pq}=\left | \begin{matrix}
a_{11}&\dots &a_{1q}&\dots &a_{1n}\\
            &\dots &           &\dots &           \\
a_{p1}&\dots &\boxed {a_{pq}}&\dots &a_{1n}\\
            &\dots &           &\dots &           \\
a_{n1}&\dots &a_{nq}&\dots &a_{nn}
\end{matrix} \right |.
$$
\begin{example}

1) For a matrix $A=(a_{ij}), i,j =1,2$
there exist four quasideterminants if the corresponding entries are invertible
$$\begin{matrix}
|A|_{11} = a_{11} - a_{12}\cdot a^{-1}_{22}\cdot a_{21},\quad
&|A|_{12}=a_{12}-a_{11}\cdot a^{-1}_{21}\cdot a_{22},\\
|A|_{21} = a_{21} - a_{22}\cdot a^{-1}_{12}\cdot a_{11},\quad &
|A|_{22}=a_{22}-a_{21}\cdot a^{-1}_{11}\cdot a_{12}.\end{matrix}
$$

2) For a matrix $A=(a_{ij}), i,j=1,2,3$ there exist nine quasideterminants
but we will write here only
$$\begin{matrix}
|A|_{11}=a_{11}-a_{12}(a_{22}-a_{23}a_{33}^{-1}a_{32})^{-1}a_{21} &-a_{12}(a_{32}-a_{33}\cdot
a_{23}^{-1} a_{22})^{-1} a_{31}\\
 \qquad -a_{13}(a_{23}-a_{22}a_{32}^{-1}a_{33})^{-1}a_{21} &-a_{13}(a_{33}-a_{32}\cdot
a^{-1}_{22}a_{23})^{-1}a_{31}\end{matrix}
$$
provided all inverses are defined.

\end{example}

Quasideterminant is not a generalization of a determinant over a commutative
ring but a generalization of a ratio of two determinants.

\begin{example} If  $A$ is a matrix over a commutative ring then
$$
|A|_{pq} = (-1)^{p+q} \frac{\det A}{\det A^{pq}}.
$$

Also, if $A$ is invertible and $A^{-1}=(b_{ij})$ then
$$b_{ij}^{-1}=|A|_{ji}$$
if the element $b_{ij}$ is invertible.

\end{example}

\begin{remark}
{\rm If each $a_{ij}$ is an invertible morphism $V_j\to V_i$ in
an additive category, then the quasideterminant $|A|_{pq}$ is also a
morphism from the object $V_q$ to the object $V_p$.
}
\end{remark}

\subsection {Elementary properties of quasideterminants}
Here is a list of elementary properties of quasideterminants.

i)  The quasideterminant $|A|_{pq}$ does not depend on the permutation
of rows and columns in the matrix $A$ if the $p$-$th$ row and the $q $-$th$
column are not changed;

ii) {\it The multiplication of rows and columns.}  Let the
matrix $B$ be constructed from the matrix $A$ by multiplication of its
$i$-$th$ row by a scalar $\lambda$  from the left.  Then
$$
 |B|_{kj}=\begin{cases} \lambda |A|_{ij} \qquad&\text{ if } k = i\\
|A|_{kj} \qquad&\text{ if } k \neq i \text{ and } \lambda \text{ is
invertible.}\end{cases}
$$
Let the matrix $C$ be constructed from the matrix A by multiplication
of its $j$-$th$ column by a scalar $\mu$ from the right.  Then
$$
|C|_{i\ell}=\begin{cases} |A|_{ij} \mu \qquad&\text{ if } \ell = j\\
|A|_{i\ell} \qquad&\text{ if } \ell \neq j \text{ and } \mu \text{ is
invertible.}\end{cases}
$$

iii) {\it The addition of rows and columns.} Let the matrix
$B$ be constructed by adding to some row of the matrix $A$ its
$k$-$th$ row multiplied by a scalar $\lambda$ from the left.  Then
$$
|A|_{ij} = |B|_{ij}, \quad  i=1, \dots k-1, k+1,\dots n,
 j=1, \dots, n.
$$
Let the matrix $C$ be constructed by addition to some column of
the matrix $A$ its $\ell$-$th$ column multiplied by a scalar
$\lambda$ from the right.  Then
$$
|A|_{ij}= |C|_{ij} , \, i=1,\dots, n, j=1,\dots ,\ell -1,\ell +1,\dots n.
$$

The following {\it homological relations} play an important role in the theory.

\begin{theorem}

a) Row homological relations:
$$
-|A|_{ij} \cdot |A^{i\ell}|^{-1}_{sj} = |A|_{i\ell}\cdot
|A^{ij}|^{-1}_{s\ell}\qquad \forall s\neq i
$$

b) Column homological relations:
$$
-|A^{kj}|^{-1}_{it} \cdot |A|_{ij} = |A^{ij}|^{-1}_{kt}\cdot |A|_{kj}
\qquad \forall r\neq j
$$
\end{theorem}

\subsection {Noncommutative Sylvester formula}
The following noncommutative version of the famous Sylvester
identity found in \cite {gr1, gr2} is closely related with  the
fundamental {\it Heredity principle} (see \cite {gr3, ggrw}).

Let $A=(A_{ij}), i,j=1,\dots, n$ be a matrix over a skew field ${\mathcal F}$.
Let $k<n-1$. Suppose $k\times k$-submatrix
$A_0=(a_{ij})$,  $i\in I_0$, $j\in J_0$ is invertible.
For $p\notin I_0$, $q\notin J_0$ construct $(k+1)\times
(k+1)$-submatrix $A_{pq}=(a_{ij})$ where $i\in I_0\cup \{p\}$,
$j\in J_0\cup \{q\}$. Set
$$
b_{pq}=|A_{pq}|_{pq}
$$
and construct matrix $B=(b_{pq})$, $p\notin I_0$, $q\notin J_0$.

We call the submatrix $A_0$ a {\it pivot} for matrix $B$.

\begin{theorem} For $s\notin I_0$, $t\notin J_0$
$$
|A|_{st} = |B|_{st}.
$$
\end{theorem}

A particular case of the theorem when $I_0=J_0=\{2,\dots , n-1\}$ is
called {\it noncommutative Lewis Carroll identity}.

\begin{example} Let $n=3$, $I_0=J_0=\{2\}$. Then $|A|_{11}$ equals to
$$
\left |\begin{matrix} \boxed {a_{11}}&a_{12}\\a_{21}&a_{22}\end{matrix}
\right | - \left |\begin{matrix} a_{12}&\boxed
{a_{13}}\\a_{22}&a_{23}\end{matrix} \right | \left |\begin{matrix}
a_{21}&a_{22}\\ \boxed {a_{31}}&a_{32}\end{matrix} \right | \left|
\begin{matrix} a_{22}&a_{23}\\ a_{32}&\boxed {a_{33}}\end{matrix} \right |.
$$

\end{example}
\subsection{Quasi-Pl\"ucker coordinates and Gauss $LDU$-factorization}
Here we remind some definitions and results from \cite {gr3, ggrw}.

Let $A=(a_{pq})$, $p=1, \dots, k$, $q=1, \dots, n$, $k< n$ be a
matrix over a skew field ${\mathcal F}$. Fix

$1\leq i, j, i_1,\dots, i_{k-1}\leq n$ such
that $i\notin I=\{ i_1,\dots, i_{k-1}\}$.  For $1\leq s\leq k$ set
$$
q^{I}_{ij}(A) =\vmatrix &a_{1i}&a_{1i_1}&\dots & a_{1i_{k-1}}\\
&{}  &\dots & {}\\
&a_{ki}&a_{ki_1} & \dots &a_{ki_{k-1}}\endvmatrix^{-1}_{si}
\cdot
\vmatrix &a_{1j}& a_{1i_1} &\dots &a_{1i_{k-1}}\\
& {}  &\dots  &{}\\
&a_{kj}& a_{ki_1} &\dots &a_{ki_{k-1}}\endvmatrix_{sj}.
$$

\begin{proposition}

$i)\ q_{ij}^I (A)$ does not depend on $s$;

$ii)\  q^I_{ij}(gA)=q^I_{ij}(A)$ for any invertible $k\times k$ matrix $g$
over ${\mathcal F}$.
\end{proposition}

We call $q^I_{ij}(A)$ left
quasi-Pl\"ucker coordinates of the matrix $A$.

In the commutative case $q^I_{ij} = \frac{p_{jI}}{p_{iI}}$, where
$p_{\alpha_1\dots \alpha_k}$ is the standard Pl\"ucker
coordinate.

Similarly, one can introduce
right quasi-Pl\"ucker coordinates.
 Consider a matrix $B=(b_{ij}), i=1,\dots, n;
j=1,\dots, k, k< n$ over a skew field ${\mathcal F}$.  Fix $1\leq i, j,
i_1,\dots, i_{k-1} \leq n$ such that $j\notin I = (i_1,\dots,
i_{k-1})$. Also fix $1\leq t\leq k$ and set
$$
r^I_{ij}(B) = \vmatrix &b_{i1} &\dots & b_{ik}\\
&b_{i_11} & \dots & b_{i_1k}\\
&{} &\dots &{}\\
&b_{i_{k-1}1} &\dots &b_{i_{k-1}k}\endvmatrix_{it} \cdot
\vmatrix &b_{j1} & \dots &b_{jk}\\
&b_{i_11} &\dots &b_{i_1 k}\\
&{} &\dots &{} \\
&b_{i_{k-1} 1} &\dots &b_{i_{k-1}k}\endvmatrix^{-1}_{jt} .
$$

\begin{proposition}

i) $r^I_{ij}(B)$ does not depend of
$t$;

ii) $r^I_{ij}(Bg) = r^I_{ij}(B)$ for any invertible
$k\times k$-matrix
$ g$ over ${\mathcal F}$.
\end{proposition}

We call $r^I_{ij}(B)$ right quasi-Pl\"ucker
coordinates of the matrix $B$.

To describe the Gauss decomposition we need the following
notations. Let $A=(a_{ij}),i,j=1,\dots,n$. Set $A^k=(a_{ij}), i,j
=k,\dots n$, $B^k=(a_{ij}), i=1,\dots n,\ j=k, \dots n$, and
$C^k=(a_{ij}), i=k,\dots n, j=1,\dots n$. These are submatrices
of sizes $(n-k+1)\times (n-k+1)$, $n\times (n-k+1)$, and
$(n-k+1)\times n$ respectively.

Suppose that the quasideterminants
$$
y_k=|A^k|_{kk},\  k=1,\dots, n
$$
are defined and invertible.
\begin{theorem}
$$
A=
\left ( \begin{matrix}1&{}&0\\
{}&\ddots&{}\\
x_{\beta \alpha}&{}&1\end{matrix}
\right )
\left ( \begin{matrix} y_1&{}&0\\
{}&\ddots&{}\\
0&{}&y_n\end{matrix} \right )
\left (\begin{matrix} 1 &{}&z_{\alpha \beta}\\
{}&\ddots&{}\\
0&{}&1\end{matrix} \right ),
$$
where
$$
\aligned
x_{\beta \alpha } &=r^{1\dots \alpha -1}_{\beta \alpha }(B^{\alpha}),
\ 1\leq \alpha <\beta\leq n\\
z_{\alpha \beta } &=q^{1\dots \alpha -1}_{\alpha\beta} ( C^{\alpha}),
\ 1\leq \alpha <\beta\leq n\\
\endaligned
$$
\end{theorem}

A noncommutative analog of the Bruhat decomposition was given in \cite {ggrw}.

\begin{example} For $n=2$
$$
A=
\left ( \begin{matrix}1&0\\
a_{21}a_{11}^{-1}&1\end{matrix}
\right )
\left ( \begin{matrix} a_{11}&0\\
0&|A|_{22}\end{matrix} \right )
\left (\begin{matrix} 1 &a_{11}^{-1}a_{12}\\
0&1\end{matrix} \right ).
$$

\end{example}

\subsection{Positive quasiminors}
\label{subsect:positive quasiminors}
Most of results in this subsection are new.

For a given matrix $x\in Mat_n(R)$ and $I,J\subset [1,n]=\{1,2,\ldots,n\}$
denote by $x_{I,J}$ the sub-matrix with the rows $I$ and the columns $J$.
And, if $|I|=|J|$, i.e., when $x_{I,J}$ is a square matrix, for any $i\in
I$, $j\in J$ denote by
$|x_{I,J}|_{i,j}$ the quasideterminant of the submatrix $x_{I,J}$ with
the marked position $(i,j)$.

Let us denote by $\Delta^i(x)$ the principal $i\times i$-quasiminor of $x\in Mat_n(R)$, i.e.,
$$\Delta^i(x)=|x_{\{1,2,\ldots,i\},\{1,2,\ldots,i\}}|_{i,i}\ .$$

The following fact is obvious.

\begin{lemma}
\label{le:permutations}
For any $I,J\subset \{1,2,\ldots,n\}$ such that $|I|=|J|=k$ and any $i\in I$, $j\in J$ there exist permutations  $u,v$ of $\{1,2,\ldots,n\}$ such that $I=u\{1,\ldots,k\}$, $J=v\{1,\ldots,k\}$, $i=u(k)$, $j=v(k)$, and for any $x\in Mat_n(R)$ we have:
$$\Delta^k(u^{-1}\cdot x\cdot v)=|x_{I,J}|_{i,j} \ .$$
(where we identified the permutations $u$ and $v$ with the corresponding $n\times n$ matrices).
\end{lemma}

\begin{definition}
{\rm For for $I,J\subset [1,n]$, $|I|=|J|$, $i\in I$, $j\in J$ define the \emph{positive quasiminor}
$\Delta^{i,j}_{I,J}$  as follows.

$$\Delta^{i,j}_{I,J}(x)=(-1)^{d_i(I)+d_j(J)}    |x_{I,J}|_{i,j}$$
where $d_i(I)$ (resp. $d_j(J)$) is the number of those elements of $I$
(resp. of $J$) which are greater than $i$ (resp. than $j$).
}

\end{definition}

The definition is motivated by the fact that for a commutative ring $R$
one has
$$\Delta^{i,j}_{I,J}(x)=\frac{\det(x_{I,J})}{\det(x_{I',J'})}\ ,$$
where $I'=I\setminus \{i\}$, $J'=J\setminus \{j\}$.
That is, a positive quasiminor is a positive ratio of minors.

Let $S_n$ be the group of permutations on $\{1,2,\ldots, n\}$.
Clearly, for any subsets $I,J\subset [1,n]$  with $|I|=|J|=k$ and elements
$i\in I$, $j\in J$ there exists a pair of
permutations $u,v\in
S_n$ such that $I=u(\{1,2,\ldots,k\})$, $J=v(\{1,2,\ldots,k\})$, $i=u(k)$, $j=v(k)$. For any such pair $u,v\in S_n$ we denote
\begin{equation}
\label{eq:positive quasiminor Sn}
\Delta^k_{u , v }:=\Delta^{i,j}_{I,J}
\end{equation}

Denote by $D_n=D_n(R)$ the set of all diagonal $n\times n$ matrices over $R$.

Clearly, positive quasiminors satisfy the relations:
\begin{equation}
\label{eq:minor properties Cartan}
\Delta^i_{u, v} (h x h') = h_{u(i)}
\Delta^i_{u, v} (x )h'_{v(i)}
\end{equation}
for $h=diag(h_1,\ldots,h_n),h'=diag(h'_1,\ldots,h'_n)\in D_n$ and
\begin{equation}
\label{eq:minor properties transposed}
\Delta^i_{u,v} (x) =  \Delta^i_{v, u} (x^T) \ ,
\end{equation}
where $x \mapsto x^T$ is the ``transpose" involutive antiautomorphism of
$Mat_n(R)$.

Let $\sigma$ be an involutive automorphism of of $Mat_n(R)$ defined by
\begin{equation}
\label{eq:tau-wnot}
\sigma(x)_{ij} = x_{n+1-i,n+1-j} \ ,
\end{equation}

The following fact is obvious.

Let  $\wnot =(n,n-1,\ldots,1)$ be the longest permutation in $S_n$.

\begin{lemma}  For any $u,v\in S_n$, and $x\in Mat_n(R)$ we have
\begin{equation}
\label{eq:minors-tau}
\Delta^i_{u, v} (\sigma(x)) =\Delta^i_{\wnot  u, \wnot  v} (x)
\end{equation}
\end{lemma}

Now we present some less obvious identities for positive
quasiminors. For each permutation $v\in S_n$ denote by $\ell(v)$
the number of inversions of $v$. Also for $i=1,2,\ldots,n-1$
denote by $s_i$ the simple transposition $(i,i+1)\in S_n$.

\begin{proposition}
\label{pr:minors-Dodgson}
Let $u,v \in S_n$ and $i \in [1,n-1]$
be such that $\l (us_i) = \l (u) + 1$ and $\l (vs_i) = \l (v) + 1$.
Then
\begin{eqnarray}
\begin{array}{l}
\label{eq:minors-Dodgson}
 \Delta^i_{us_i , v s_i }
= \Delta^i_{us_i , v }(\Delta^i_{u, v})^{-1}\Delta^i_{u , vs_i }
+  \Delta^{i+1}_{u , v  } \ ,
\end{array}
\end{eqnarray}
$$
(\Delta^i_{u s_i, v} )^{-1}\Delta^{i+1}_{u , v}
=(\Delta^i_{u, v} )^{-1}\Delta^{i+1}_{u s_i , v}
\ , \Delta^{i+1}_{u , v} (\Delta^i_{u , vs_i} )^{-1}
=\Delta^{i+1}_{u , vs_i }(\Delta^i_{u, v} )^{-1},$$
$$
 \Delta^{i+1}_{u , v} (\Delta^{i+1}_{u s_i , v})^{-1}
=\Delta^i_{u s_i, v}(\Delta^i_{u, v} )^{-1}
\ , (\Delta^{i+1}_{u , vs_i })^{-1}\Delta^{i+1}_{u , v}
=(\Delta^i_{u, v} )^{-1}\Delta^i_{u , vs_i}. $$
\end{proposition}

\begin{proof} Clearly, the fourth and the fifth identities follow from the second and the third. Using Lemma \ref{le:permutations} and the Gauss factorization it suffices to take $u=v=1$, $i=1$ in the first three identities, i.e., work with $2\times 2$ matrices.
Then the first three identities will take respectively the following obvious form:
$$x_{22}=x_{21}x_{11}^{-1}x_{12}+\left|\begin{matrix}
x_{11}&x_{12}\\
x_{21}&\boxed {x_{22}}~
\end{matrix}\right |\ ,$$

$$
x_{21}^{-1}\left|
\begin{matrix}
x_{11}&x_{12}\\
x_{21}&\boxed {x_{22}}~
\end{matrix}\right |
=-x_{11}^{-1}\left|
\begin{matrix}
x_{11}&\boxed {x_{12}}~\\
x_{21}&x_{22}
\end{matrix}
\right |
\ ,
\left|
\begin{matrix}
x_{11}&x_{12}~\\
x_{21}&\boxed {x_{22}}
\end{matrix}\right |x_{12}^{-1}
=-\left|
\begin{matrix}
x_{11}&x_{12}\\
\boxed {x_{21}}&x_{22}
\end{matrix}
\right |x_{11}^{-1} \ . $$

\end{proof}

One can prove the next proposition presenting some generalized
Pl\"ucker relations.

\begin{proposition}
\label{pr:minors-Plucker}
Let $u, v \in S_n$ and $i \in [1,n-2]$.
If $\l (us_i s_{i+1} s_i) = \l (u) + 3$, then
\[
 \Delta^{i+1}_{us_{i+1} , v} =
 \Delta^{i+1}_{us_i s_{i+1} , v} +
  \Delta^i_{us_{i+1} s_i , v }(\Delta^i_{u s_i, v })^{-1}\Delta^{i+1}_{u , v }\ .
\]

If $\l (vs_i s_{i+1} s_i) = \l (v) + 3$, then
\[
 \Delta^{i+1}_{u , vs_{i+1}} =
 \Delta^{i+1}_{u , vs_i s_{i+1}} +
 \Delta^{i+1}_{u , v }(\Delta^i_{u , vs_i })^{-1} \Delta^i_{u , vs_{i+1} s_i }\ .
\]
\end{proposition}

\section{Basic factorizations in  $GL_n({\mathcal F})$}

\label{sec:basic factorizations}

For $i,j=1,2,\ldots,n$ denote by $E_{ij}$ the $n\times n$ matrix unit in the
intersection of the $i$-th row and the $j$-th column.

Then we abbreviate $E_i:=E_{i,i+1}$ for $i=1,\ldots,n-1$.

The matrix units $E_1,\ldots,E_{n-1}$ satisfy the relations:
$E_i^2=0$ for $i=1,\ldots,n-1$ and
$$E_iE_j=E_jE_i$$ if $|i-j|\ge 2$,
$$E_iE_{i\pm 1}E_i=0 \ .$$

Let $\ii=(i_1,\ldots,i_m)$ be a sequence of indices $i_k\in \{1,2,\ldots,n-1\}$
and $x=(x_{ij})$, $i,j=1,\dots , n$ be an $n\times n$-matrix over a skew
field ${\mathcal
F}$.
For such an $\ii$ and $x$ let us write
the formal factorization,
\begin{equation}
\label{eq:X factored}x=(1+t_1 E_{i_1})(1+t_2 E_{i_2})\cdots (1+t_m
E_{i_m}) \ ,
\end{equation}
where all $t_k$ belong to the skew field ${\mathcal F}$.

\medskip

Let $k_{ij}$ can be the position of $i$-th occurrence of the index $j-i$ in the  sequence
$\ii=(1,\ldots,n-1;1,\ldots,n-2;\ldots;1,2;1)$. That is,
$$k_{ij}=n(i-1)-\binom{i+1}{2}+j$$
for $1\le i< j\le n$.

\begin{proposition} Let $\ii=(1,\ldots,n-1;1,\ldots,n-2;\ldots;1,2;1)$.
We set temporarily $t_{ij}:=t_{k_{ij}}$ for $1\le i< j\le n$
(where $t_k$ are as in the factorization (\ref{eq:X factored})).
Then the matrix entries of the product $x$ satisfy ($1\le i\le n-k\le
n-1$):
\begin{equation}
\label{eq:xik}
x_{i,i+k}=\sum\limits_{1\le i_1\le i_2\le \cdots \le i_k\le n+1-i-k} t_{i_1,i_1+i}t_{i_2,i_2+i+1}t_{i_3,i_3+i+2}\cdots
t_{i_k,i_k+i+k-1}  \ .
\end{equation}

\end{proposition}

\begin{remark}
{\rm
In particular, after the specialization $t_{k_{ij}}:=y_{j}$ (for $1\le i< j\le n$)
in (\ref{eq:xik}) for some elements $y_2,\ldots,y_n$, we obtain:
$$x_{i,i+k}=\sum\limits_{i<  j_1<j_2<\cdots <j_k\le n} y_{j_1}y_{j_2}\cdots
y_{j_k} \ .
$$
That is, each matrix entry of so specialized matrix $x$ is an elementary
symmetric function in $y_2,\ldots,y_n$.
}
\end{remark}

\begin{proposition} The system (\ref{eq:xik}) has a unique solution of the form:
\begin{equation}
\label{eq:quasi-ratio}
t_{ij}=|x_{i,j-1}|_{j-i, n-i+1}\cdot |x_{ij}|_{j-i+1, n-i+1}^{-1}
\end{equation}
for $1\le i<j\le n$, where $x_{ij}$ is the $i\times i$-submatrix of $x$
with the
rows $\{j-i+1,\ldots,j\}$ and the columns $\{n-i+1,\ldots,n\}$.

\end{proposition}

\begin{proof} First of all, we have the relations
$$x_{i,n}=t_{1,i+1}t_{1,i+2}\cdots t_{1,n}$$
for $i=1,\ldots,n-1$. Therefore,
$$t_{1,i+1}=x_{i,n}x_{i+1,n}^{-1}$$
for all $j=2,\ldots,n$ which is verifies (\ref{eq:quasi-ratio}).

Let us define a sequence $x^{(0)}, x^{(1)}, \ldots, x^{(n-1)}$ of matrices
inductively by setting $x^{(0)}=I$ and
$$x^{(m)}=(I+t_{n-m,n+1-m}E_{12})(I+t_{n-m,n+2-m}E_{23})\cdots
(I+t_{n-m,n}E_{m,m+1})\cdot x^{(m-1)}$$
for $m=1,2,\ldots,n-1$.

Clearly, $x^{(n-1)}=x$.

\begin{lemma} One has for all $i\le j\le m+1\le n$:
\begin{equation}
\label{eq:quasi-coefficients}
x^{(m)}_{ij}=|x_{ij}^m|_{ij} \ .
\end{equation}
where $x_{ij}^m$ is the $(n-m)\times (n-m)$ submatrix of
$x$ with the rows $\{i,i+1,\ldots,i+n-m-1\}$ and the columns
$\{j;m+2,m+3,\ldots,n \}$.
\end{lemma}

\begin{proof} We proceed by induction on $n-m$. By definition of
$x^{(m)}$,
we have a recursion for the matrix entries of $x^{(m)}$:
$$x_{i,j}^{(m)}= t_{n-m,i+n-m}\cdot x_{i+1,j}^{(m)}+x^{(m-1)}_{ij}$$
for $1\le i\le j\le m+1$.

Taking $j=m+1$, we obtain
\begin{equation}
\label{eq:t extreme}
t_{n-m,i+n-m}=x_{i,m+1}^{(m)}\cdot (x_{i+1,m+1}^{(m)})^{-1}
\end{equation}
Therefore,
$$x^{(m-1)}_{ij}=x_{i,j}^{(m)}-x_{i,m+1}^{(m)}\cdot (x_{i+1,m+1}^{(m)})^{-1}x_{i+1,j}^{(m)}
=\left |\begin{matrix}\boxed {x_{i,j}^{(m)}}&x_{i,m+1}^{(m)}\\
x_{i+1,j}^{(m)}&x_{i+1,m+1}^{(m)}
\end{matrix}\right | \ . $$

Furthermore, let us use the inductive hypotheses precisely in the form
(\ref{eq:quasi-coefficients}).
Then, by the above,
$$x^{(m-1)}_{ij}=\left |\begin{matrix}\boxed
{|x_{ij}^m|_{ij}}&|x_{i,m+1}^m|_{i,
m+1}\\
|x_{i+1,j}^m|_{i+1, j}&|x_{i+1,m+1}^m|_{i+1, m+1} \end{matrix}\right
|\ . $$

Using the Sylvester formula (Theorem 1.6)  with $A=x_{i,j,m-1}$ and $A_0$
being a
submatrix of $x$ with the rows $\{i+1,\ldots,i+n-m-1\}$ and the columns
$\{m+2,m+3,\ldots,n \}$, we obtain:
$$x^{(m-1)}_{ij}=
\left |\begin{matrix}\boxed {|x_{ij}^m|_{i,j}}&|x_{i,m+1}^m|_{i,m+1}\\
|x_{i+1,j}^m|_{i+1,j}&|x_{i+1,m+1}^m|_{i+1,m+1} \end{matrix}\right|
=|x_{ij}^{m-1}|_{ij} \ . $$

This finishes the induction. The lemma is proved. \end{proof}

Finally, using (\ref{eq:quasi-coefficients}), (\ref{eq:t extreme}), and
the fact that $x_{i,m+1}=x_{i,m+1}^m$ for $m=0,1,\ldots,n-1$, we obtain
(\ref{eq:quasi-ratio}).

The proposition is proved.
\end{proof}

Another natural factorization of generic matrices is given by the
following theorem.

For a generic matrix $x=(x_{ij})$, $i,j=1,\dots , n$
over a skew field ${\mathcal F}$ define the sequence
of rational functions $t_{m,k}=t_{m,k}(x)$, $1\le m\le k\le n-1$
by the formula:

$$
t_{m,k}=\left | \begin{matrix}
x_{1, k-m+1}&\dots  & x_{1k}\\
                    &\dots  &            \\
x_{m, k-m+1}&\dots  &\boxed {x_{mk}}
\end{matrix} \right  |^{-1}
\cdot
\left | \begin{matrix}
x_{1, k-m+2}&\dots  & x_{1, k+1}\\
                    &\dots  &            \\
x_{m, k-m+2}&\dots  & \boxed {x_{m, k+1}}
\end{matrix} \right  |$$
Clearly, in terms of positive quasiminors, one has:
$$t_{m,k}=(\Delta^{m,k}_{\{1,\ldots,m\},\{k-m+1,\ldots,k\}})^{-1}\Delta^{m,k+1}_{\{1,\ldots,m\},\{k-m+2,\ldots,k+1\}} \ .$$

Then define a sequence of matrices $x(m,k)=(x_{ij}^{(m,k)})$, $1\le m\le k\le n-1$  by the inductive formula:
$$x(1,n-1)=x\cdot (1-t_{1, n-1}E_{n-1})$$
$$x(m,k)=x(m, k+1)\cdot (1-t_{m, k}E_k )$$
$$x(m+1, n-1)=x(m, m)\cdot (1-t_{m+1, n-1}E_{n-1}) \ .$$
In other words,
$$x(m,k)\cdot  \prod_{(i,j)\preceq (m,k)} (1+t_{i, j}E_j)=x \ ,$$
where the order $\prec$ on all pairs $(m,k)$, $1\le m\le k\le n-1$, is defined by: $(i,j)\prec (m,k)$ if and only if either $i<m$ or $i=m$, $j>k$.

\begin{theorem}
\label{th:upper factorization}
(a) For a generic matrix $x=(x_{ij})$, $i,j=1,\dots , n$
over a skew field ${\mathcal F}$ one has
$$x_{ij}^{(m,k)}=0$$ for all $i,j$ such that $(i,j-1) \preceq (m,k)$
(i.e., for $i<j, i<m$ and for $i=m, j>k$). In particular, $x(n-1,n-1)$ is
lower triangular.

b) The  entries $x_{ij}^{(m,k)}$ are given by the following formulas: For $i\geq m$, $2\leq j\leq k$
$$
x_{ij}^{(m,k)}=\left | \begin{matrix}
x_{1, j-m+1}&\dots  & x_{1j}\\
                    &\dots  &            \\
x_{m-1, j-m+1}&\dots  & x_{m-1,j}\\
x_{i, j-m+1}&\dots  &\boxed {x_{i,j}}
\end{matrix} \right  |,
$$
for $i> m$, $j>k$
$$
x_{ij}^{(m,k)}=\left | \begin{matrix}
x_{1, j-m}&\dots  & x_{1,j}\\
                    &\dots  &            \\
x_{m, j-m}&\dots  & x_{m,j}\\
x_{i, j-m}&\dots  &\boxed {x_{i,j}}
\end{matrix} \right  |,
$$
and $x_{ij}^{(m,k)}=x_{ij}$ otherwise.

\end{theorem}

\begin{proof} It is enough to show that matrices $x(m, k)$ satisfy
conditions i)-iv) listed below.

i)  $x_{ij}^{(m,k)}=0$ for $i<j, i<m$ and for $i=m, j>k$,

ii) $x(1,n-1)=x(1+E_{n-1}t_{1, n-1})$,

iii) $x(m,k)=x(m, k+1))(1+E_kt_{m, k})$,

iv)  $x(m+1, n-1)=x(m, m)(1+E_{n-1}t_{m+1, n-1})$.

We proceed by induction over a totally ordered set of indices
$(1,n-1),\dots , (1,1)$, $(2, n-1),\dots , (2,2)$, $\dots $,
$(n-1, n-1)$.

It is easy to check that the entries of matrix $x(1,n-1)$ satisfy
conditions i)-iv). Suppose that these conditions are satisfied for for matrix $x(m, l)$.
We consider then two cases: $l>m$ and $l=m$.

If $l>m$ then $l=k+1$ for $k\geq m)$. Define matrix $x(m,k)$ by formula iii).
Evidently, the corresponding entries of matrices $x(m,k)$ and  $x(m,k+1)$
coincide except the entries with indices $i,k$ for $i\geq m$ which are given by
the formula
$$x_{ik}^{(m,k)}=x_{ik}^{(m,k+1)}t_{m,k}+x_{ik+1}^{(m,k+1)}.$$

For $i\geq m$ the product $x_{ik}^{(m,k+1)}t_{m,k}$ equals to
$$
- \left | \begin{matrix}
x_{1, k-m+1}&\dots  & x_{1,j}\\
                    &\dots  &            \\
x_{m-1, k-m+1}&\dots  & x_{m-1,j}\\
x_{i, k-m+1}&\dots  &\boxed {x_{i,j}}
\end{matrix} \right  |
\left | \begin{matrix}
x_{1, k-m+1}&\dots  & x_{1k}\\
                    &\dots  &            \\
x_{m, k-m+1}&\dots  &\boxed {x_{mk}}
\end{matrix} \right  |^{-1}
\left | \begin{matrix}
x_{1, k-m+2}&\dots  & x_{1, k+1}\\
                    &\dots  &            \\
x_{m, k-m+2}&\dots  & \boxed {x_{m, k+1}}
\end{matrix} \right  |.
$$
According to the homological relations for quasideterminants
the last expression can be written as
$$
-\left | \begin{matrix}
x_{1, k-m+1}&\dots  & x_{1,j}\\
                    &\dots  &            \\
x_{m-1, k-m+1}&\dots  & x_{m-1,j}\\
\boxed {x_{i, k-m+1}}&\dots  &x_{i,j}
\end{matrix} \right  |
\left | \begin{matrix}
x_{1, k-m+1}&\dots  & x_{1k}\\
                    &\dots  &            \\
\boxed {x_{m, k-m+1}}&\dots  &x_{mk}
\end{matrix} \right  |^{-1}
\left | \begin{matrix}
x_{1, k-m+2}&\dots  & x_{1, k+1}\\
                    &\dots  &            \\
x_{m, k-m+2}&\dots  & \boxed {x_{m, k+1}}
\end{matrix} \right  |.
$$
It follows that the element $x_{ik}^{(m,k)}=x_{ik}^{(m,k+1)}t(m,k)+x_{ik+1}^{(m,k+1)}$
is zero for $i=m$. If $i>m$
$$
x_{ij}^{(m,k)}=\left | \begin{matrix}
x_{1, k-m+1}&\dots  & x_{1,k+1}\\
                    &\dots  &            \\
x_{m, k-m+1}&\dots  & x_{m,k+1}\\
x_{i, k-m+1}&\dots  &\boxed {x_{i,k+1}}
\end{matrix} \right  |.
$$
It follows from  the Sylvester identity applied to the corresponding matrix with the
pivot equal to
$$
\left ( \begin{matrix}
x_{1, k-m+2}&\dots  & x_{1,k}\\
                    &\dots  &            \\
x_{m-1, k-m+2}&\dots  &x_{m-1,k}
\end{matrix} \right  ).
$$

It shows that the entries of matrix $x(m,k)$ satisfy part b) of the theorem.

If $l=m$ one can check in a similar way that the entries of matrix $x(m+1, n-1)$
satisfy  part b) of the theorem.

The theorem is proved.
\end{proof}

\begin{remark}
{\rm  It follows from the proof that matrices $x(m,k)$
and elements $t_{m,k}$ are uniquely defined.
}
\end{remark}

\begin{example} Let $n=3$. Then $$t_{1,2}=-x_{12}^{-1}x_{13},\ \ \
 t_{1,1}=-x_{11}^{-1}x_{12},$$
$$ t_{2,2}=-\left | \begin{matrix}
x_{11}&x_{12}\\
x_{21}&\boxed {x_{22}}
\end{matrix} \right |^{-1}
\cdot
\left | \begin{matrix}
x_{12}&x_{13}\\
x_{22}&\boxed {x_{23}}
\end{matrix} \right |,
$$
$$
x(1,2)=\left ( \begin{matrix}
x_{11}&x_{12}&0\\
x_{21}&x_{22}&
\left | \begin{matrix}
x_{12}&x_{13}\\
x_{22}&\boxed {x_{23}}
\end{matrix} \right |\\
x_{31}&x_{32}&
\left | \begin{matrix}
x_{12}&x_{13}\\
x_{32}&\boxed {x_{33}}
\end{matrix} \right |
\end{matrix} \right )
$$
$$
x(1,1)=\left ( \begin{matrix}
x_{11}&0&0\\
x_{21}&
\left | \begin{matrix}
x_{11}&x_{12}\\
x_{21}&\boxed {x_{22}}
\end{matrix} \right |
&
\left | \begin{matrix}
x_{12}&x_{13}\\
x_{22}&\boxed {x_{23}}
\end{matrix} \right |\\
x_{31}&
\left | \begin{matrix}
x_{11}&x_{12}\\
x_{31}&\boxed {x_{32}}
\end{matrix} \right |
&
\left | \begin{matrix}
x_{12}&x_{13}\\
x_{32}&\boxed {x_{33}}
\end{matrix} \right |
\end{matrix} \right )
$$
$$
x(2,2)=\left ( \begin{matrix}
x_{11}&0&0\\
x_{21}&
\left | \begin{matrix}
x_{11}&x_{12}\\
x_{21}&\boxed {x_{22}}
\end{matrix} \right |
&0\\
x_{31}&
\left | \begin{matrix}
x_{11}&x_{12}\\
x_{31}&\boxed {x_{32}}
\end{matrix} \right |
&|x|_{33}
\end{matrix} \right )
$$
\end{example}


\section{Examples}
\label{sec:examples}

\subsection{A factorization in the Borel subgroup of $GL_3({\mathcal F})$}

\label{subsect:SL3}

Let us write the formal factorization
$$\begin{pmatrix}
x_{11}& x_{12}& x_{13}\\
0&x_{22}&x_{23}\\
0&0&x_{33}\\
\end{pmatrix}=
\begin{pmatrix}
x_{11}& 0& 0\\
0&x_{22}&0\\
0&0&x_{33}\\
\end{pmatrix}
\begin{pmatrix}
1&t_{12}+t_{23}&t_{12}t_{13}\\
0&1&t_{13}\\
0&0&1
\end{pmatrix}$$
$$=\begin{pmatrix}
x_{11}& 0& 0\\
0&x_{22}&0\\
0&0&x_{33}\\
\end{pmatrix}
\begin{pmatrix}
1&t_{12}&0\\
0&1&0\\
0&0&1\\
\end{pmatrix}
\begin{pmatrix}
1&0&0\\
0&1&t_{13}\\
0&0&1\\
\end{pmatrix}
\begin{pmatrix}
1&t_{23}&0\\
0&1&0\\
0&0&1\\
\end{pmatrix}
$$

(assuming that all $x_{ij}, t_{ij}$ are elements of a skew field ${\mathcal F}$).

Then we can express $t_{ij}$ as follows.
$$t_{13}=x_{22}^{-1}x_{23},
~t_{12}=x_{11}^{-1}x_{13}x_{23}^{-1}x_{22},
~t_{23}=x_{11}^{-1}x_{12}-x_{11}^{-1}x_{13}x_{23}^{-1}x_{22}=x_{11}^{-1}\left|
\begin{matrix}
\boxed {x_{12}}&x_{13}\\
x_{22}&x_{23}
\end{matrix}
\right |.$$

\begin{remark}
{\rm The above factorization exists (and, therefore, is unique)
if and only if each of $x_{11},x_{22},x_{33}$, and $x_{23}$ is invertible.
}
\end{remark}

\subsection{A factorization in $GL_3({\mathcal F})$}

\label{subsect:GL3}

Let us write the formal factorization over a skew field ${\mathcal F}$.
$$x=\begin{pmatrix}
x_{11}& x_{12}& x_{13}\\
x_{21} &x_{22}&x_{23}\\
x_{31} &x_{32} &x_{33}\\
\end{pmatrix}=hx_{-2}(t_1)x_{-1}(t_2)x_{-2}(t_3)x_2(t_4)x_1(t_5)x_2(t_6)$$
where
$$h=\begin{pmatrix}
h_1& 0& 0\\
0&h_2&0\\
0&0&h_3\\
\end{pmatrix},
x_1(t)=\begin{pmatrix}
1&t&0\\
0&1&0\\
0&0&1
\end{pmatrix},
x_2(t)=\begin{pmatrix}
1&0&0\\
0&1&t\\
0&0&1
\end{pmatrix},
$$
$$x_{-1}(t)=\begin{pmatrix}
t^{-1}&0&0\\
1&t&0\\
0&0&1
\end{pmatrix},
x_{-2}(t)=\begin{pmatrix}
1&0&0\\
0&t^{-1}&0\\
0&1&t
\end{pmatrix}.
$$

Then we can express $h_i$ and $t_k$ as follows.

$$h_3=x_{31}, h_2=-\left|
\begin{matrix}
x_{21}&\boxed {x_{22}}\\
x_{31}&x_{32}
\end{matrix}
\right |, h_1=\left |
\begin{matrix}
x_{11}&x_{12}&\boxed {x_{13}}\\
x_{21}&x_{22}&x_{23}\\
x_{31}& x_{32}&x_{33}
\end{matrix}
\right |,$$

$$t_6=x_{12}^{-1}x_{13}, t_5=x_{11}^{-1}x_{12},$$
$$t_4=(x_{22}-x_{21}x_{11}^{-1}x_{12})^{-1}(x_{23}-x_{22}x_{12}^{-1}x_{13})=\left|
\begin{matrix}
x_{11}&x_{12}\\
x_{21}&\boxed {x_{22}}
\end{matrix}
\right |^{-1}\left|
\begin{matrix}
x_{12}&x_{13}\\
x_{22}&\boxed {x_{23}}
\end{matrix}
\right |,$$
$$t_1=-x_{21}^{-1}\left|
\begin{matrix}
x_{21}&\boxed {x_{22}}\\
x_{31}&x_{32}
\end{matrix}
\right |,
t_2= x_{11}^{-1}\left |
\begin{matrix}
x_{11}&x_{12}&\boxed {x_{13}}\\
x_{21}&x_{22}&x_{23}\\
x_{31}& x_{32}&x_{33}
\end{matrix}
\right |,
t_3=-\left |
\begin{matrix}
x_{11}&\boxed {x_{12}}\\
x_{21}&x_{22}
\end{matrix} \right |^{-1}
\left |\begin{matrix}
x_{11}&x_{12}&\boxed {x_{13}}\\
x_{21}&x_{22}&x_{23}\\
x_{31}& x_{32}&x_{33}
\end{matrix}
\right |.$$

In fact, if we define a sequence of matrices
$$x^{(5)}=x\cdot x_2(t_6)^{-1}, x^{(4)}=x^{(5)}x_1(t_5)^{-1},
x^{(3)}=x^{(4)}x_1(t_4)^{-1} \ $$
then $x^{(k-1)}$  will have exactly one more zero entry in the upper part
than $x^{(k)}$:
$$x^{(5)}=\begin{pmatrix}
x_{11}& x_{12}& 0\\
x_{21} &x_{22}&x'_{23}\\
x_{31} &x_{32} &x'_{33}\\
\end{pmatrix}, x^{(4)}=\begin{pmatrix}
x_{11}& 0& 0\\
x_{21} &x''_{22}&x'_{23}\\
x_{31} &x''_{32} &x'_{33}\\
\end{pmatrix}, x^{(3)}=\begin{pmatrix}
x_{11}& 0& 0\\
x_{21} &x''_{22}& 0\\
x_{31} &x''_{32} &x'''_{33}\\
\end{pmatrix}.$$
This determines $t_6,t_5,t_4$.

And the rest of parameters $h_1,h_2,h_3, t_1,t_2, t_3$ are obtained from the equation:
$$x^{(3)}=hx_{-2}(t_1)x_{-1}(t_2)x_{-2}(t_3)=\begin{pmatrix}
h_1t_2^{-1}& 0& 0\\
h_2t_1^{-1} &h_2t_1^{-1}t_2t_3^{-1}& 0\\
h_3 &h_3(t_1+t_2t_3^{-1}) &h_3t_1t_3\\
\end{pmatrix}$$

\subsection{A factorization in the unipotent subgroup of $GL_4({\mathcal F})$}

\label{subsect:SL4}

\label{sect:Factorization SL4}
Let us write the formal factorization
$$\begin{pmatrix}
1& x_{12}& x_{13}&x_{14}\\
0&1&x_{23}&x_{24}\\
0&0&1&x_{34}\\
0&0&0&1
\end{pmatrix}=
\begin{pmatrix}
1&t_{12}+t_{23}+t_{34}&t_{12}t_{13}+t_{12}t_{24}+t_{23}t_{24}&t_{12}t_{13}t_{14}\\
0&1&t_{13}+t_{24}&t_{13}t_{14}\\
0&0&1&t_{14}\\
0&0&0&1\end{pmatrix}$$
$$=\begin{pmatrix}
1&t_{12}&0&0\\
0&1&0&0\\
0&0&1&0\\
0&0&0&1
\end{pmatrix}
\begin{pmatrix}
1&0&0&0\\
0&1&t_{13}&0\\
0&0&1&0\\
0&0&0&1
\end{pmatrix}
\begin{pmatrix}
1&0&0&0\\
0&1&0&0\\
0&0&1&t_{14}\\
0&0&0&1\end{pmatrix}
\times
$$

$$\times \begin{pmatrix}
1&t_{23}&0&0\\
0&1&0&0\\
0&0&1&0\\
0&0&0&1\end{pmatrix}\begin{pmatrix}
1&0&0&0\\
0&1&t_{24}&0\\
0&0&1&0\\
0&0&0&1\end{pmatrix}
\cdot \begin{pmatrix}
1&t_{34}&0&0\\
0&1&0&0\\
0&0&1&0\\
0&0&0&1\end{pmatrix}
 $$
(assuming that all $x_{ij}, t_{ij}$ are elements of a  skew field ${\mathcal F}$).

Then we can express $t_k$ as follows.
$$t_{14}=x_{34},
~t_{13}=x_{24}x_{34}^{-1},~t_{12}=x_{14}x_{24}^{-1},
~t_{24}=x_{23}-x_{24}x_{34}^{-1}=\left|
\begin{matrix}
\boxed {x_{23}}&x_{24}\\
1&x_{34}
\end{matrix}
\right |,$$
$$t_{23}=(x_{13}-x_{14}x_{24}^{-1}
x_{23})(x_{23}-x_{24}x_{34}^{-1})^{-1}=
\left |
\begin{matrix}
\boxed {x_{13}}&x_{14}\\
x_{23}&x_{24}
\end{matrix}
\right | \left |\begin{matrix}\boxed {x_{23}}&x_{24}\\  1&x_{34}\end{matrix}\right
|^{-1},~$$
$$t_{34}=x_{12}-x_{13}(x_{23}-x_{24}x_{34}^{-1})^{-1}+x_{14}x_{34}^{-1}
(x_{23}-x_{24}x_{34}^{-1})^{-1}=\left |
\begin{matrix}
\boxed {x_{12}}&x_{13}&x_{14}\\
1&x_{23}&x_{24}\\
0&1&x_{34}
\end{matrix}
\right |.$$

\begin{remark}
{\rm The above factorization exists (and, therefore, is unique)
if and only if $x_{24}$, $x_{34}$, and
$\left |\begin{matrix}\boxed {x_{23}}&x_{24}\\  1&x_{34}\end{matrix}\right |$
are invertible in ${\mathcal F}$.
}
\end{remark}

\section{Double Bruhat cells in $GL_n({\mathcal F})$ and their factorizations}
\label{sec:double bruhat cells}

\subsection{Structure of $GL_n({\mathcal F})$}
\label{sec:geometric prelims}
Throughout this and the next section we denote $G:=GL_n({\mathcal F})$ and will use the abbreviation (for $a,b\in \ZZ$):

\begin{equation}
\label{eq:ab}
[a,b]=\begin{cases} \{a,a+1,\ldots,b\} & \text{if $a\le b$} \\
\emptyset & \text{otherwise }
\end{cases}
\end{equation}
Let $U$ (resp. $U^-$) be the upper (resp. lower) unitriangular subgroup of $G$.
For $i \in [1,r]$, we define the elementary unitriangular matrices  $x_i (t)$ and $y_i (t)$ by:
$$x_i (t) = I+tE_i \,\, , y_i (t) = I+tF_i$$
for $i\in [1,n-1]$, where $E_i=E_{i,i+1}$, $F_i=E_{i+1,i}$ are the corresponding matrix units (in the notation of Section \ref{sec:basic factorizations}).

Let $H$ denote the subgroup of all diagonal matrices in $G$.
Let $B$ (resp. $B^-$) be the subgroup of all upper (resp. lower) triangular matrices in $G$.
Clearly,  $B = HU$, $B^- = HU^-$, and  $H=B^-\cap B$.

We denote by $G_0=B^-U$ the open subset of elements $x\in G$ that
have Gaussian $LDU$-decomposition; this (unique) decomposition will be written as
$x = [x]_-  [x]_+ \,$ (where $[x]_-\in B^-$, but not necessarily in $U^-$)
For any $x$ in the Gauss cell
$G_0=B^-\cdot U$ denote by $[x]_0$ the diagonal component of the Gauss $LDU$-factorization.
In particular, $[x]_0=[[x]_-]_0$ for any $x\in G_0$.

For each $i\in [1,n-1]$, let $\varphi_i: GL_2({\mathcal F}) \to G$ denote
the  embedding corresponding to the $2\times 2$ block at the intersection of the $i$-th and $(i+1)$st rows and the $i$-th and $(i+1)$st columns.
Thus we have
$$x_i (t) = \varphi_i \mat{1}{t}{0}{1}, \,\, y_i (t) = \varphi_i \mat{1}{0}{t}{1} \ .$$
We also set
$$h_i(t) = \varphi_i \mat{t}{0}{0}{t^{-1}} \in H, x_{-i}(t)=\varphi_i \mat{t^{-1}}{0}{1}{t}$$
for any $i$ and any $t \in {\mathcal F}^\times$. By definition,
$$x_{-i}(t)=y_i(t)h_i(t^{-1}) =h_i(t^{-1})y_i(t^{-1})\ .$$
More generally, it is  easy to see that for each $i\in [1,n-1]$ and any diagonal matrix $h=diag(h_1,\ldots,h_n)\in H$ one has:
\begin{equation}
\label{eq:Borel relation}
hx_i(t)h^{-1}=x_i(h_ith_{i+1}^{-1}), h^{-1}y_i(t)h=y_i(h_{i+1}^{-1}th_i)
\end{equation}
Hence
\begin{equation}
\label{eq:Cartan relations}
h_j(s)x_i(t)=x_i(s^{\varepsilon_{ji}}ts^{\varepsilon_{ij}})h_j(s),\\
y_i(t)h_j(s)=h_j(s)y_i(s^{\varepsilon_{ij}}ts^{\varepsilon_{ji}})
\end{equation}
for any $i,j\in [1,n-1]$, where $\varepsilon_{ij}=\delta_{ij}-\delta_{i,j-1}$.


\begin{lemma}
\label{le:x-i, xi} $~$

(i) For each $i\in [1,n-1]$ we have: $x_{-i}(s)x_i(t)=x_i(s^{-1}t(s+t)^{-1})x_{-i}(s+t)$.

(ii) For each $i\in [1,n-2]$ we have:
$x_{-i}(s)x_{i+1}(t)=x_{i+1}(st)x_{-i}(s)$.

(iii) For each $i\in [2,n-1]$ we have:
$x_{-i}(s)x_{i-1}(t)=x_{i-1}(ts)x_{-i}(s)$.

(iv) For any $i, j\in [1,n-1]$ such that $|i-j|>1$ we have:
$$x_{-i}(s)x_j(t)=x_j(t)x_{-i}(s) \ .$$
\end{lemma}

\begin{proof}
Part (i) follows from the obvious identity:
$$
\begin{pmatrix}
s^{-1}  &   0 \\
1   & s
\end{pmatrix}
\begin{pmatrix}
1  &    t \\
0   & 1
\end{pmatrix}
=
\begin{pmatrix}
s^{-1}  &   s^{-1}t \\
1   & s+t
\end{pmatrix}=
\begin{pmatrix}
1  &    s^{-1}t(s+t)^{-1} \\
0   & 1
\end{pmatrix}\begin{pmatrix}
(s+t)^{-1}  &   0 \\
1   & s+t
\end{pmatrix}$$
for $s,t\in {\mathcal F}^\times$.

Part (ii) follows from
$$
\begin{pmatrix}
s^{-1}  &   0 & 0\\
1   & s & 0 \\
0 & 0 & 1
\end{pmatrix}
\begin{pmatrix}
1& 0 & 0 \\
0 & 1  &    t \\
0& 0   & 1
\end{pmatrix}
=
\begin{pmatrix}
s^{-1}  &   0 & 0\\
1   & s & st \\
0 & 0 & 1
\end{pmatrix}
=
\begin{pmatrix}
1& 0 & 0 \\
0 & 1  &    st \\
0& 0   & 1
\end{pmatrix}
\begin{pmatrix}
s^{-1}  &   0 & 0\\
1   & s & 0 \\
0 & 0 & 1
\end{pmatrix}$$
for $s,t\in {\mathcal F}^\times$.

Part (iii) follows from
$$
\begin{pmatrix}
1 & 0 & 0 \\
0 & s^{-1}  &   0 \\
0 & 1   & s
\end{pmatrix}
\begin{pmatrix}
1  & t & 0\\
0  & 1 & 0\\
0  & 0 & 1
\end{pmatrix}
=
\begin{pmatrix}
1 & t & 0 \\
0 & s^{-1}  &   0 \\
0 & 1   & s
\end{pmatrix}
=
\begin{pmatrix}
1  & ts & 0\\
0  & 1 & 0\\
0  & 0 & 1
\end{pmatrix}
\begin{pmatrix}
1 & 0 & 0 \\
0 & s^{-1}  &   0 \\
0 & 1   & s
\end{pmatrix}$$
for $s,t\in {\mathcal F}^\times$.

And part (iv) is obvious.
\end{proof}

The symmetric group $S_n$ of $G$ is naturally embedded into $G$ via
$$(i,i+1)\mapsto \varphi_i\mat{0}{1}{1}{0}$$
for $i\in [1,n-1]$. We also define a representative
$\overline {s_i}$ of the transposition $(i,i+1)$ by
\[
\overline {s_i} = \varphi_i \mat{0}{-1}{1}{0} \, .
\]
The elements $\overline {s_i}$ satisfy the braid relations in~$W$;
thus the representative $\overline w$ can be unambiguously defined for any
$w \in W$ by requiring that
$\overline {uv} = \overline {u} \cdot \overline {v}$
whenever $\l (uv) = \l (u) + \l (v)$.

\subsection{Bruhat cells and Double Bruhat cells}
\label{sec:reduced Bruhat}

The group $G$ has two \emph{Bruhat decompositions},
with respect to opposite Borel subgroups $B$ and $B^-\,$:
$$G = \bigcup_{u \in S_n} B u B = \bigcup_{v \in S_n} B^- v B^-  \ . $$

Now define the {\it Schubert cell} $U(w):=wU^-w^{-1}\cap U$ for $w\in S_n$. Then the following obvious fact demonstrates that the {\it Bruhat cells} $B u B$ and $B^- v B^-$ behave similarly to their commutative counterparts.

\begin{lemma}
\label{le:Bruhat}
(a) For each $u\in S_n$ one has:
$$BuB=U(u)uB=BuU(u),~U\overline uU=U(u)\overline uU=UuU(u^{-1})\ .$$
(b) For each $v\in S_n$ one has:
$$B^-vB^-=B^-U(v)\overline v=B^-U(v){\overline {v^{-1}}}^{\,-1}=\overline v U(v^{-1})B^-={\overline {v^{-1}}}^{\,-1}U(v^{-1})B^- \ .$$

\end{lemma}

\begin{definition}
{\rm
For any permutations $u,v\in S_n$ define the \emph{double Bruhat cell}~$G^{u,v}$ by
$G^{u,v} = B u B  \cap B^- v B^- \,$.
}
\end{definition}

In this section  we shall concentrate on the following subset $L^{u,v} \subset G^{u,v}$
which we call a \emph{reduced double Bruhat cell}:
\begin{equation}
\label{eq:reduced cell}
L^{u,v} = U \overline u U  \cap B^- v B^- \ .
\end{equation}


\begin{remark}
{\rm In the commutative case the reduced double Bruhat cells are simplectic leafs of the Poisson-Lie structure on $GL_n(\CC)$ (see e.g., \cite{kogzel}). These cells also emerge in the study of {\it total positivity} (\cite{bz}) on $GL_n$.
}
\end{remark}

The equations defining $L^{u,v}$ inside $G^{u,v}$ look as follows.

\begin{proposition}
\label{pr:Luv equations}
An element $x \in G^{u,v}$ belongs to $L^{u,v}$ if and only if
$[{\overline u}^{\ -1} x]_0 = 1$, or equivalently if $\Delta^i_{u , e} (x) = 1$
for each $i \in  [1,n]$.
\end{proposition}

The maximal torus $H$ acts freely on $G^{u,v}$ by left (or right) translations,
and $L^{u,v}$ is a section of this action.
Thus $L^{u,v}$ is naturally identified with $G^{u,v}/H$, and all properties of $G^{u,v}$
can be translated in a straightforward way into the corresponding properties of $L^{u,v}$.

A \emph{double reduced word} for a pair $u,v \in S_n$
is a reduced word for an element $(u,v)$
of the  group $S_n \times S_n$.
To avoid confusion, we will use the indices
$- 1, \ldots, - r$ for the simple reflections
in the first copy of $W$, and $1, \ldots, r$ for the second copy.
A double reduced word for $(u,v)$ is simply a shuffle of a reduced
word $\ii$ for $u$ written in the alphabet
$[- 1, - r]$ (we will denote such a word by $- \ii$) and a reduced word $\ii'$
for $v$ written in the alphabet $[1,r]$.
We denote the set of double reduced words for $(u,v)$ by $R(u,v)$.

For any sequence $\ii= (i_1, \ldots, i_m)$ of indices
from the alphabet $[1,r] \cup [- 1, - r]$,
let us define the \emph{product map} $x_\ii: ({\mathcal F}^\times)^m \to G$ by
\begin{equation}
\label{eq:productmap}
x_\ii (t_1, \ldots, t_m) =  x_{i_1} (t_1) \cdots x_{i_m} (t_m) \, .
\end{equation}

\subsection{Factorization problem for reduced double Bruhat cells}
\label{sec:factorization}
In this section, we address the following
\emph{factorization problem} for $L^{u,v}$: for any double reduced word
$\ii \in R(u,v)$, find explicit formulas for the inverse birational isomorphism
$x_\ii^{-1}$ between $L^{u,v}$ and $({\mathcal F}^\times)^m$,
thus expressing the \emph{factorization parameters} $t_k$
in terms of the product~$x = x_\ii (t_1, \ldots, t_m) \in L^{u,v}$.

\begin{definition}
{\rm
Let $x \mapsto x^\iota$ be the involutive antiautomorphism of $G$
given by
$$x^\iota=J_nx^{-1}J_n $$
for any $x\in G$, where
$J_n=diag(-1,1,-1,\ldots,(-1)^{n})$.
}

\end{definition}

We will refer to the anti-automorphism $x \mapsto x^\iota$ as to the {\it positive inverse} in $G$.
It is easy to see that
\begin{equation}
\label{eq:iota}
a^\iota = a^{-1} \quad (a \in H) \ , \quad x_i (t)^\iota = x_i (t) \ ,
\quad y_{i} (t)^\iota = y_{i} (t) \ .
\end{equation}

The following fact is a direct noncommutative analogue of
Theorem 1.6 from \cite {fz}.

\begin{lemma}
\label{le:iota Bruhat}For any $u,v\in S_n$ one has:
$$(BuB)^\iota=Bu^{-1}B,~(U\overline uU)^\iota=U\overline {u^{-1}}U,~ (B^-vB^-)^\iota=B^-v^{-1}B\ . $$
In particular, $(G^{u,v})^\iota=G^{u^{-1},v^{-1}}$.
\end{lemma}

\begin{definition}
\label{def:twist}
{\rm
For any $u,v\in W$, the twist map
$\psi^{u,v}:L^{u,v}\to G$ is defined by
\begin{equation}
\label{eq:psi-u,v-x prime}
\psi^{u,v}(x) =
([x \overline{v^{-1}}\,]_-)^\iota \, (x^\iota)^{-1} \,
([\overline{u}^{\ -1}x]_+)^{\iota} \ .
\end{equation}

}

\end{definition}

\begin{theorem}
\label{th:psi-regularity}
The twist map $\psi^{u,v}$ is an isomorphism
between $L^{u,v}$ and  $L^{v, u}$.
The inverse isomorphism is $\psi^{v,u}$.
\end{theorem}

\begin{proof} The proof essentially follows the pattern of the commutative case
from \cite{fz} and \cite{bz}. We  need the following obvious fact.
\begin{lemma}
\label{le:twist}
The twist map $\psi^{u,v}$ is satisfies:
\begin{equation}
\label{eq:psi-u,v-x}
\psi^{u,v}(x) =
[(\overline{v} x^{\iota})^{-1}]_+ \, \overline{v} \,
([\overline{u}^{\ -1} x]_+)^{\iota} = ([x \overline{v^{-1}}\,]_-)^\iota \, \overline{u^{-1}}^{\,-1} \,
[\overline{u^{\ -1}} ((x)^\iota)^{-1}]_-\ .
\end{equation}
The restriction of $\psi^{u,v}$ to $L^{u,v}\cap B^- U$ is a
map $L^{u,v}\cap B^- U\to L^{u,v}\cap B^-U$ given by the formula:
\begin{equation}
\label{eq:psi-u,v-x G0}
\psi^{u,v}(x_-\cdot x_+) =
([x_+ \overline{v^{-1}} ]_-)^\iota \cdot
([\overline{u}^{\ -1} x_-]_+)^{\iota} \ .
\end{equation}

In particular, the twist map $\psi^{u,e}:L^{u,e}\to L^{e,u}$ is given by
\begin{equation}
\label{eq:psi wnot-e}
\psi^{u,e}(x) = ([\overline{u}^{\ -1} x]_+)^{\iota} \ .
\end{equation}
And $\psi^{e,v}:L^{e,v}\to L^{v,e}$ is given by
\begin{equation}
\label{eq:psi e-wnot}
\psi^{e,v}(x) = ([x\overline{v^{-1}}]_- )^\iota\ .
\end{equation}
\end{lemma}

The formula (\ref{eq:psi-u,v-x}) guarantees that
$\psi^{u,v}(L^{u,v})\subset U\overline vU\cap B^-uB^-=L^{v,u}$, i.e., $\psi^{u,v}$ is
a well-defined map $L^{u,v}\to L^{v,u}$.

Finally, we prove that $\psi^{v,u}$ is the inverse of $\psi^{u,v}$, i.e.,
$\psi^{v,u}\circ \psi^{u,v}=id$. Given $x\in L^{u,v}$, denote $y=\psi^{u,v}(x)$.
By definition (\ref{eq:psi-u,v-x prime}), we have

$$y =
([x \overline{v^{-1}}\,]_-)^\iota \, (x^\iota)^{-1} \,
([\overline{u}^{\ -1}x]_+)^{\iota} \ .$$
Or, equivalently,
$$(y^\iota)^{-1}=(([x \overline{v^{-1}}\,]_-)^{-1} \, x \,
([\overline{u}^{\ -1}x]_+)^{-1} \ , $$
and
$$x=[x \overline{v^{-1}}\,]_- \, (y^\iota)^{-1} \,
[\overline{u}^{\ -1}x]_+ \ . $$

Since
$$
\psi^{v,u}(y) =
([y \overline{u^{-1}}\,]_-)^\iota \, (y^\iota)^{-1} \,
([\overline{v}^{\ -1}y]_+)^{\iota} \ ,
$$
in order to prove that $\psi^{v,u}(y)=x$ it suffices to show that
$$([y \overline{u^{-1}}\,]_-)^\iota=[x \overline{v^{-1}}\,]_-, ~([\overline{v}^{\ -1}y]_+)^{\iota}=
[\overline{u}^{\ -1}x]_+ \ ,$$
or, equivalently,
\begin{equation}
\label{eq: x y plus minus}
[y \overline{u^{-1}}\,]_-=([x \overline{v^{-1}}\,]_-)^\iota, ~[\overline{v}^{\ -1}y]_+=[\overline{u}^{\ -1}x]_+)^\iota \ .
\end{equation}
Let us prove the first identity (\ref {eq: x y plus minus}).
Denote temporarily $z=([x \overline{v^{-1}}\,]_-)^\iota$. Then (\ref{eq:psi-u,v-x}) implies that
$$y \overline{u^{-1}}=z \cdot \overline{u^{-1}}^{\,-1} \,
[\overline{u^{\ -1}} ((x)^\iota)^{-1}]_-\overline{u^{-1}} \ .$$
According to Lemma \ref{le:iota Bruhat}, for any $x\in U\overline u U$
we have: $((x)^\iota)^{-1}\in U \overline{u^{-1}}^{\,-1} U$, and, furthermore, by
Lemma \ref{le:Bruhat}(a),
$\overline{u^{\ -1}}x^\iota\in \overline{u^{\ -1}}U(u) \overline{u^{-1}}^{\,-1} U \subset U^-\cdot U$,
and $[\overline{u^{\ -1}} ((x)^\iota)^{-1}]_-\in \overline{u^{\ -1}}U(u) \overline{u^{-1}}^{\,-1}$.
Hence $\overline{u^{-1}}^{\,-1} \,[\overline{u^{\ -1}} ((x)^\iota)^{-1}]_-\overline{u^{-1}}\in U$.
Therefore,
$$[y \overline{u^{-1}}]_-=[z \cdot \overline{u^{-1}}^{\,-1} \,
[\overline{u^{\ -1}} ((x)^\iota)^{-1}]_-\overline{u^{-1}}]_-=[z]_-=z \ .$$
This proves the first identity in (\ref {eq: x y plus minus}).
Now let us prove the second identity in (\ref {eq: x y plus minus}).
Again denote temporarily $t=([\overline{u}^{\ -1}x]_+)^\iota$. Then (\ref{eq:psi-u,v-x}) implies that
$$\overline{v}^{\ -1}y=\overline{v}^{\ -1}[(\overline{v} x^{\iota})^{-1}]_+ \, \overline{v} \,
\cdot t \ .$$
According to Lemma \ref{le:Bruhat}(b), for any $x\in B^-v B^-$ one has
$x \overline {v^{-1}}\in B\cdot U(v)$, and $[(\overline{v} x^{\iota})^{-1}]_+\in U(v)$.
Hence $\overline{v}^{\ -1}[(\overline{v} x^{\iota})^{-1}]_+ \, \overline{v}\in \overline{v}^{\ -1}U(v) \overline{v}\subset U^-$.
Therefore,
$$[\overline{v}^{\ -1}y]_+=[\overline{v}^{\ -1}[(\overline{v} x^{\iota})^{-1}]_+ \, \overline{v} \,
\cdot t]_+=[t]_+=t \ .$$
This proves the second identity in (\ref {eq: x y plus minus}).

Theorem \ref{th:psi-regularity} is proved.  \end{proof}

Now let us fix a pair $(u,v) \in S_n \times S_n$ and a double reduced word
$\ii = (i_1, \ldots, i_m)\in R(u,v)$.
Recall that $\ii$ is a shuffle of a reduced word for $u$ written in the alphabet
$[- 1, - r]$ and a reduced word for $v$ written in the alphabet $[1,r]$.
In particular, the length $m$ of $\ii$ is equal to $\l (u) + \l (v)$.

We will use the convention that $s_{-i}=1$ for each $i\in [1,n-1]$.
For $k\in [1,m]$, denote
\begin{equation}
\label{eq:ulessk}
u_{\geq k} =
s_{-i_m}s_{-i_{m-1}}\cdots s_{-i_k} \ , \quad u_{> k} = s_{-i_m}s_{-i_{m-1}}\cdots s_{-i_{k+1}} \ ,
\end{equation}
\begin{equation}
\label{eq:vlessk}
v_{\leq k} = s_{i_1}s_{i_2}\cdots s_{i_k} \ , \quad v_{<k} = s_{i_1}s_{i_2}\cdots s_{i_{k-1}} \ .
\end{equation}
For example, if
$\ii= (-2, 1, -3, 3, 2, -1, -2, 1, -1)$,
then, say, $u_{\geq 7} = s_1 s_2 \,$ and $v_{<7} = s_1 s_3 s_2\,$.

Now we are ready to state our solution to the factorization problem.

\begin{theorem}
\label{th:t-through-x}
Let $\ii = (i_1, \ldots, i_m)$ be
a double reduced word for $(u,v)$, and suppose an element
$x \in L^{u,v}$ can be factored as
$x= x_{i_1} (t_1) \cdots x_{i_m} (t_m)$,
with all $t_k$ nonzero elements of ${\mathcal F}$.
Then the factorization parameters $t_k$
are determined by the following formula:
\begin{equation}
\label{eq:t-through-x}
t_k =\begin{cases}
\Delta^i_{v_{<k} , u_{> k}}(y)^{-1}\Delta^i_{v_{<k} , u_{\geq k} } (y)=\Delta^{i+1}_{v_{<k}, u_{\geq k}}(y)^{-1}\Delta^{i+1}_{v_{<k} , u_{> k} } (y) & \text{if $i_k<0$}\\
\Delta^i_{v_{\le k}, u_{> k}} (y)^{-1}\Delta^{i+1}_{v_{<k} , u_{> k}} (y)
=\Delta^i_{v_{< k}, u_{> k}} (y)^{-1}\Delta^{i+1}_{v_{\le k} , u_{> k}} (y) & \text{if $i_k>0$}
\end{cases}
\end{equation}
where $y=\psi^{u,v} (x)$ and $i=|i_k|$.
\end{theorem}

\begin{proof} First, let us list some important properties of positive quasiminors. Recall
that in Section \ref{subsect:positive quasiminors}, for $i\in [1,n]$ we defined the \emph{principal quasi-minor} $\Delta^i$ by:
$$\Delta^i(x)=|x_{[1,i],[1,i]}|_{i,i}$$
for any $x\in G$, where $x_{[1,i],[1,i]}$ denotes the principal $i\times i$ submatrix of $x$. In particular, $\Delta^1(x)=x_{11}$ and $\Delta^n(x)=|x|_{n,n}$.

The following fact is obvious.

\begin{lemma} The principal quasi-minors are invariant under the left multiplication by $U^-$ and the right multiplication by $U$, i.e.,
$$\Delta^i(x_-xx_+)=\Delta^i(x)$$
for any $x_+\in U$, $x_-\in U^-$, $x\in G$ (in particular, $\Delta^i(x)=\Delta^i([x]_0)=([x]_0)_{ii}$). Furthermore, for any $u,v\in S_n$ one has
\begin{equation}
\label{eq:Delta-general}
\Delta^i_{u , v } (x) =
\Delta^i(\overline {u}^{\ -1}
   x \overline v) \ .
\end{equation}
Also one has:
$$ \Delta^i_{u, v} (x^\iota)=\Delta^{n+1-i}_{\wnot  v\wnot , \wnot  u\wnot } (x)^{-1} \ .$$

\end{lemma}

We will prove (\ref{eq:t-through-x}) by the induction in the $l(u)+l(v)$. The base of the induction with $u=v=e$ is obvious.

We will consider the following four cases:

Case I. $u\ne e$, $v\ne e$ and $\ii$ is {\it separated}, i.e, $-i_1,\ldots,-i_{\ell}\in [1,n-1]$ and $i_{\ell+1},\ldots,i_m\in [1,n-1]$ for some $\ell$, or, equivalently, $u=s_{-i_1} \cdots s_{-i_\ell}$ and $v=s_{i_{\ell+1}}\cdots s_{i_m}$.

Case II. $u\ne e$, $v\ne e$ and $\ii$ is not separated.

Case III $u=e$, $v\ne e$.

Case IV. $u\ne e$, $v=e$.

Consider Case I first.

Denote
$$x_-:=x_{i_1}(t_1)\cdots x_{i_\ell}(t_\ell), ~x_+:=x_{i_{\ell+1}}(t_{\ell+1})\cdots x_{i_m}(t_m) \ .$$ Clearly,
$x_-\in L^{u,e}$, $x_+\in L^{e,v}$, and
$x=x_-\cdot x_+\in L^{u,v}$.
Furthermore,
the inductive hypothesis (\ref{eq:t-through-x}) for $x_-$ says that:
$$t_k =\Delta^i_{e , u_{> k}}(y_+)^{-1}\Delta^i_{v_{<k} , u_{\geq k} } (y_+)=\Delta^{i+1}_{e, u_{\geq k}}(y_+)^{-1}\Delta^{i+1}_{e , u_{> k} } (y_+)$$
 for $k\in [1,\ell]$, where $y_+=\psi^{u,e}(x_-)$, $i=|i_k|$.

And the inductive hypothesis (\ref{eq:t-through-x}) for $x_+$ says that
$$t_k=\Delta^i_{v_{\le k}, e} (y_-)^{-1}\Delta^{i+1}_{v_{<k} , e} (y_-)
=\Delta^i_{v_{< k}, e} (y_-)^{-1}\Delta^{i+1}_{v_{\le k} , e} (y_-)
$$
for $k\in [\ell+1,m]$, where $y_-=\psi^{e,v}(x_+)$, $i=|i_k|$.

According to (\ref{eq:psi-u,v-x G0}), (\ref{eq:psi wnot-e}), and (\ref{eq:psi e-wnot}),
$$\psi^{u,v}(x)=([x_+ \overline{v^{-1}} ]_-)^\iota \cdot
([\overline{u}^{\ -1} x_-]_+)^{\iota}=y_-y_+ \ .$$
Note also that $\Delta^j_{e,w}(y_+)=\Delta^j_{e,w}(y_-y_+)$ and $\Delta^j_{w,e}(y_-)=\Delta^j_{w,e}(y_-y_+)$ for any $w\in S_n$ and $j\in [1,n]$.
Finally, taking into the account that $v_{\le k}=v_{<k}=e$  for each $k\le \ell$, and $u_{\ge k}=u_{>k}=e$ for each $k>\ell$, we obtain  (\ref{eq:t-through-x}) for $x=x_-x_+$. This finishes Case I.

Now consider Case II. We say that given $\ii$, a pair $(i_\ell,i_{\ell+1})$ is an inversion if $i_\ell>0$ and $i_{\ell+1}<0$. Clearly, $\ii$ has no inversions if and only if $\ii$ is separated. Here we will proceed by the induction in the number of inversions.  The base of the induction is the already considered Case I -- no inversions. Assume that $\ii'$ has an inversion $(i'_\ell,i'_{\ell+1})=(i,-j)$, where $i,j\in [1,n-1]$. Let $\ii$ be obtained form $\ii'$ by switching $i_\ell$ and $i_{\ell+1}$, that is, $\ii$ has one inversion less than $\ii$. According to the inductive hypothesis, (\ref{eq:t-through-x}) holds for the factorization (relative to $\ii$):
$$x=x_{i_1}(t_1)\cdots x_{i_{\ell-1}}(t_{\ell-1})x_{-j}(t_\ell)x_i(t_{\ell+1})x_{i_{\ell+2}}(t_{\ell+2})\cdots x_{i_m}(t_m) \ .$$

Note that, according to Lemma \ref{le:x-i, xi},
$$x_{-j}(t_\ell)x_i(t_{\ell+1})=x_i(t'_\ell)x_{-j}(t'_{\ell+1})\ ,$$
where
\begin{equation}
\label{eq:tl plus 1 prime}
(t'_\ell,t'_{\ell+1})=\begin{cases} (t_{\ell+1},t_\ell) & \text{if $|i-j|>1$} \\
(t_\ell t_{\ell+1},t_\ell) & \text{if $i-j=1$} \\
(t_{\ell+1}t_\ell,t_\ell) & \text{if $i-j=-1$} \\
(t_\ell^{-1}t_{\ell+1}(t_\ell+t_{\ell+1})^{-1},t_\ell+t_{\ell+1}) & \text{if $i=j$} \\
\end{cases}\ .
\end{equation}
We have to prove that each of the  parameters $t_1,\ldots, t_{\ell-1},t'_\ell,t'_{\ell+1},t_{l+2},\ldots,t_m$ in the factorization (relative to $\ii'$)
$$x=x_{i_1}(t_1)\cdots x_{i_{\ell-1}}(t_{\ell-1})x_{i}(t'_\ell)x_{-j}(t'_{\ell+1})x_{i_{\ell+2}}(t_{\ell+2})\cdots x_{i_m}(t_m)$$
is given by (\ref{eq:t-through-x}) for $\ii'$.

Obviously, if $k\ne \ell, \ell+1$, then $v_{<k}^{\ii'}=v_{<k}^\ii$, $v_{\le k}^{\ii'}=v_{\le k}^\ii$, $u_{>k}^{\ii'}=u_{>k}^\ii$, $u_{\ge k}^{\ii'}=u_{\ge k}^\ii$. Therefore,  each $t_k$, $k\ne \ell, \ell+1$ in the latter decomposition is given by (\ref{eq:t-through-x}) for $\ii'$. It remains prove that $t'_{\ell}$ and $t'_{\ell+1}$ are both given by (\ref{eq:t-through-x}) for $\ii'$. Denote temporarily $u'=u_{>l}$, $v'=v_{<l}$ so that
(taking into account that $i_\ell=i'_{\ell+1}-j$, $i_{\ell+1}=i'_\ell=i$) we have $v^\ii_{\le  \ell}=v^\ii_{< \ell+1}=v'$, $v^\ii_{\le  \ell+1}=v's_i$,
$u^\ii_{\ge  \ell+1}=u^\ii_{> \ell+1}=u'$, $u^\ii_{\ge  \ell}=u's_j$. Therefore, (\ref{eq:t-through-x}) for $\ii$ with $k=\ell$ and $k=\ell+1$ becomes (with the convention  $y=\psi^{u,v} (x)$, $y'=\overline {v'}^{-1}y\overline u'$):

$$t_\ell =  \Delta^j_{e , e}(y')^{-1}\Delta^j_{e , s_j } (y')=\Delta^{j+1}_{e, s_j}(y')^{-1}\Delta^{j+1}_{e , e} (y')\ ,$$
$$t_{\ell+1}=
\Delta^i_{s_i, e} (y')^{-1}\Delta^{i+1}_{e , e} (y')
=\Delta^i_{e, e} (y')^{-1}\Delta^{i+1}_{s_i , e} (y') \ .$$

Taking into the account that $v^{\ii'}_{<  \ell}=v'$,  $v^{\ii'}_{\le  \ell}=v^{\ii'}_{< \ell+1}=v^{\ii}_{\le  \ell+1}=v's_i$, $u^{\ii'}_{>  \ell+1}=u'$,
$u^{\ii'}_{\ge  \ell+1}=u^{\ii'}_{> \ell+1}=u^{\ii'}_{\ge  \ell}=u's_j$ we have only to prove that
\begin{equation}
\label{eq:tl prime}
t'_\ell=\Delta^i_{s_i, s_j} (y')^{-1}\Delta^{i+1}_{e , s_j} (y')\ ,
=\Delta^i_{e, s_j} (y')^{-1}\Delta^{i+1}_{s_i , s_j} (y') \ .
\end{equation}
\begin{equation}
\label{eq:tl plus one prime} t'_{\ell+1}=\Delta^j_{s_i , e}(y')^{-1}\Delta^j_{s_i , s_j } (y')=\Delta^{j+1}_{s_i,s_j}(y')^{-1}\Delta^{j+1}_{s_i , e } (y')
\end{equation}

Consider the following four sub-cases:

1. $|i-j|>1$. Then clearly,  $\Delta^i_{s_i, s_j} (y')=\Delta^i_{s_i, e} (y')$, $\Delta^{i+1}_{e, s_j} (y')=\Delta^{i+1}_{e, e} (y')$, and
 $\Delta^j_{s_i, s_j} (y')=\Delta^j_{e, s_j} (y')$, $\Delta^j_{s_i, e} (y')=\Delta^j_{e, e} (y')$.  Finally, by (\ref{eq:tl plus 1 prime}), $t'_{\ell}=t_\ell t_{\ell+1}$ and $t'_{\ell+1}=t'_\ell$. All these immediately imply (\ref{eq:tl prime}) and (\ref{eq:tl plus one prime}).

2. $j=i-1$. According to (\ref{eq:tl plus 1 prime}),
$$t_\ell =  \Delta^{i}_{e, s_j}(y')^{-1}\Delta^{i}_{e , e} (y'), t_{\ell+1}=\Delta^i_{e, e} (y')^{-1}\Delta^{i+1}_{s_i , e} (y') \ , $$
$$t'_\ell=t_\ell t_{\ell+1}=\Delta^{i}_{e, s_j}(y')^{-1}\Delta^{i+1}_{s_i , e} (y')=\Delta^{i}_{e, s_j}(y')^{-1}\Delta^{i+1}_{s_i , s_j} (y')\ , $$
which proves (\ref{eq:tl prime}). Similarly,  we obtain
$$t'_{\ell+1}=t_\ell =  \Delta^j_{e , e}(y')^{-1}\Delta^j_{e , s_j } (y')=\Delta^j_{s_i , e}(y')^{-1}\Delta^j_{s_i , s_j } (y')\ ,$$
which proves (\ref{eq:tl plus one prime}).

3. $j=i+1$. According to (\ref{eq:tl plus 1 prime}),
$$t'_\ell=t_{\ell+1}t_\ell=\left(\Delta^i_{s_i, e} (y')^{-1}\Delta^{i+1}_{e , e} (y')\right)\left(\Delta^j_{e , e}(y')^{-1}\Delta^j_{e , s_j } (y')\right) =\Delta^i_{s_i, e} (y')^{-1}\Delta^j_{e , s_j } (y')\ ,$$
which proves (\ref{eq:tl prime}) because $\Delta^i_{s_i, e} (y')=\Delta^i_{s_i, s_{i+1}} (y')$. Similarly, we obtain
$$t'_{\ell+1}=t_\ell = \Delta^{j+1}_{e, s_j}(y')^{-1}\Delta^{j+1}_{e , e} (y')=\Delta^{j+1}_{s_i, s_j}(y')^{-1}\Delta^{j+1}_{s_i , e} (y')\ ,$$
which proves (\ref{eq:tl plus one prime}).

4. $i=j$. According to (\ref{eq:tl plus 1 prime}),
$$t'_{\ell+1}=t_\ell+t_{\ell+1}=\Delta^i_{e , e}(y')^{-1}\Delta^i_{e , s_i } (y')+
\Delta^i_{s_i, e} (y')^{-1}\Delta^{i+1}_{e , e} (y')=$$
$$\Delta^i_{s_i, e} (y')^{-1}\left(\Delta^i_{s_i, e} (y')\Delta^i_{e , e}(y')^{-1}\Delta^i_{e , s_i } (y')+
\Delta^{i+1}_{e , e} (y')\right)=\Delta^i_{s_i, e} (y')^{-1}\Delta^i_{s_i, s_i} (y')$$
by (\ref{eq:minors-Dodgson}). This proves (\ref{eq:tl plus one prime}).

Furthermore, according to (\ref{eq:tl plus 1 prime}),
$$t'_lt'_{l+1}=t_\ell^{-1}t_{\ell+1}=\left( \Delta^i_{e , s_i}(y')^{-1}\Delta^i_{e , e } (y')\right)\left(\Delta^i_{e, e} (y')^{-1}\Delta^{i+1}_{s_i , e} (y') \right)=\Delta^i_{e , s_i}(y')^{-1}\Delta^{i+1}_{s_i , e} (y') \ .$$
Therefore, using already proved (\ref{eq:tl plus one prime}), we obtain:
$$t'_\ell=\Delta^i_{e , s_i}(y')^{-1}\Delta^{i+1}_{s_i , e} (y')(t'_\ell)^{-1}=\Delta^i_{e , s_i}(y')^{-1}\Delta^{i+1}_{s_i , e} (y')\left(\Delta^{i+1}_{s_i,e}(y')^{-1}\Delta^{i+1}_{s_i , s_i } (y')\right)$$
$$=\Delta^i_{e , s_i}(y')^{-1}\Delta^{i+1}_{s_i , s_i } (y') \ ,$$
which proves (\ref{eq:tl prime}). This finishes Case II.

\smallskip

Now we consider Case III: $\ii=(i_1,\ldots,i_m)$, where all $i_k>0$, i.e, $\ii$ is a reduced word for $v$. And let $i=i_m$ so that $v=v's_i$ and $l(v)=l(v')+1$. Let
$$x=x_i(t_1)\cdots x_{i_m}(t_m), ~x'=x_i(t_1)\cdots x_{i_{m-1}}(t_{m-1})x_{-i}(t_m^{-1}) \ .$$
It is easy to see that
$$x\overline {s_i}x_i(t_m^{-1})=x' \ .$$
Indeed, this follows from
\begin{equation}
\label{eq:x plus minus si}
x_{-i}(t^{-1})=x_i(t)\overline s_ix_i(-t^{-1}) \ ,
\end{equation}
which, in its turn,  follows  from the obvious identity:
$$
\begin{pmatrix}
1  &    t \\
0   & 1
\end{pmatrix}
\begin{pmatrix}
0  &   -1 \\
1   & 0
\end{pmatrix}=
\begin{pmatrix}
t  &   -1 \\
1   & 0
\end{pmatrix}=
\begin{pmatrix}
t  &   0 \\
1   & t^{-1}
\end{pmatrix}
\begin{pmatrix}
1  &    -t^{-1} \\
0   & 1
\end{pmatrix} \ .$$

Note that $x'$ is factored along the reduced word
$\ii'=(i_1,\ldots,i_{m-1};-i)$ for $(s_i,v')$. Therefore, we can
use the already proved Case II for  the $\ii'$-factorization of $x'$.
Formula (\ref{eq:t-through-x}) for the factorization parameters
$t_1,\ldots,t_{m-1},t_m^{-1}$ of $x'$ takes the form:
$$t_k =\Delta^{i_k}_{v_{\le k}, s_i} (y')^{-1}\Delta^{i_k+1}_{v_{<k} , s_i} (y')=\Delta^{i_k}_{v_{< k}, s_i} (y')^{-1}\Delta^{i_k+1}_{v_{\le k} , s_i} (y')$$
for $k\in [1,m-1]$, and
$$t_m^{-1}=\Delta^i_{v' , e}(y')^{-1}\Delta^i_{v' , s_i } (y')=\Delta^{i+1}_{v', s_i}(y)^{-1}\Delta^{i+1}_{v' , e} (y')\ ,$$
where $y'=\psi^{s_i,v'} (x')$.

Clearly, in order to finish Case III, i.e., to verify formula (\ref{eq:t-through-x}) for the $\ii$-factorization parameters $t_1,\ldots,t_m$ of $x\in L^{e,v}$, it will suffice to prove that for any $w\in S_n$, $j\in [1,n]$ one has:
$$\Delta^{j}_{w, s_i}(y')=\Delta^j_{w, e}(y)\ ,$$
where $y=\psi^{e,v} (x)$. Note that $\Delta^{j}_{w, s_i}(y')=\Delta^{j}_{w, e}(y'\overline{s_i})$. Thus, it will suffice to prove
$$[y'\overline{s_i}]_-=y \ .$$

Taking into the account that
$$x\overline {s_i}x_i(t^{-1})=x' \ ,$$
all we need to prove is the following fact.

\begin{lemma}  Let   $v=v's_i$ for some $i$ such that $l(v)=l(v')+1$. Then for any $x''\in L^{e,v'}$ and any $t\in {\mathcal F}^\times$ one has
$$\psi^{e,v} (x''x_i(t))=[\psi^{s_i,v'} (x''x_{-i}(t^{-1}))\overline{s_i}]_-\ .$$
\end{lemma}

\begin{proof} Indeed, by Lemma \ref{le:twist},
$$\psi^{e,v} (x''x_i(t))=([x''x_i(t) \overline{v^{-1}}\,]_-)^\iota \ .$$
Using (\ref{eq:x plus minus si}), we obtain:
$$x''x_i(t) \overline{v^{-1}}=x''x_i(t)\overline s_i\overline{{v'}^{-1}}=x''x_{-i}(t^{-1})x_i(-t^{-1})\overline{{v'}^{-1}}=x'x_{-i}(t^{-1})\overline{{v'}^{-1}}u_+$$
for some $u_+\in U$ .

Therefore,
$$[x''x_i(t)\overline{v^{-1}}]_-=[x''x_{-i}(t^{-1})\overline{{v'}^{-1}}u_+]_-
=[x''x_{-i}(t^{-1})\overline{{v'}^{-1}}]_-\ .$$
Summarizing, we obtain:

$$\psi^{e,v} (x''x_i(t))=([x''x_{-i}(t^{-1}) \overline{{v'}^{-1}}\,]_-)^\iota $$
On the other hand, by the second identity of (\ref{eq:psi-u,v-x}) we have for any $x'\in L^{s_i,v'}$:
$$[\psi^{s_i,v'}(x')\overline s_i]_- =[([x' \overline{v^{'-1}}\,]_-)^\iota \, \overline{s_i}^{\,-1} \,
[\overline{s_i} ((x')^\iota)^{-1}]_+\overline {s_i}]_-=([x' \overline{v^{'-1}}\,]_-)^\iota$$
because $z=[\overline{s_i} ((x')^\iota)^{-1}]_+\in U\cap \varphi_i(GL_2)$ and, therefore,  $\overline{s_i}^{\,-1} \,
z\overline {s_i}\in B^-$.
Thus, taking $x'=x''x_{-i}(t^{-1})$, we obtain
$$[\psi^{s_i,v'}(x''x_{-i}(t^{-1})\overline s_i]_- =([x''x_{-i}(t^{-1}) \overline{v^{'-1}}\,]_-)^\iota=\psi^{e,v} (x'x_i(t)) \ .$$

Lemma is proved. \end{proof}

This finishes Case III.

Case IV is almost identical to the Case III.

Therefore, Theorem \ref{th:t-through-x} is proved.
\end{proof}

\begin{remark}
{\rm The commutative version of (\ref{eq:t-through-x}) is
\begin{equation}
t_k=\begin{cases}
\displaystyle
{\frac{\Delta_{v_{<k} \omega_{i}, u_{\geq k} \omega_{i}} (y)}{\Delta_{v_{< k} \omega_{i}, u_{> k}
\omega_{i}}(y)}} &\text{if $i_k < 0$}\\
& \\
\displaystyle
{\frac{\Delta_{v_{<k} \omega_{i-1}, u_{\geq k} \omega_{i-1}} (y)\Delta_{v_{<k}
\omega_{i+1}, u_{\geq k} \omega_{i+1}} (y)}{\Delta_{v_{<k}
\omega_{i}, u_{\geq k} \omega_{i}} (y) \Delta_{v_{\le k} \omega_{i}, u_{> k} \omega_{i}}(y)}} &\text{if $i_k > 0$}
\end{cases}
\end{equation}

}
\end{remark}

\subsection{Factorizations of $G^{u,v}$}
In this section we extend the result of Theorem \ref{th:t-through-x} to factorizations in $G^{u,v}$. In order to do so we first have to  extend the twist $\psi^{u,v}$ to an isomorphism $G^{u,v}\widetilde \to G^{v,u}$
(which we will denote in the same way) by
\begin{equation}
\label{eq:psi-u,v-x general}
\psi^{u,v}(hx)=h\psi^{u,v}(x)
\end{equation}
for any $h\in H$ and any $x\in L^{u,v}$.

In fact, formula (\ref{eq:psi-u,v-x general}) means that the twist $\psi^{u,v}$
is a {\it left $H$-equivariant} map $G^{u,v}\widetilde \to G^{v,u}$.

Recall that for any $g$ in the Gauss cell
$G_0=B^-\cdot U$ we denote by $[g]_0$ the diagonal component of the Gauss factorization.

\begin{lemma}
\label{le:twist general} The general twist $\psi^{u,v}:G^{u,v}\widetilde \to G^{v,u}$ is given by:
\begin{equation}
\label{eq:twist general}
\psi^{u,v}(g) =
u([\overline{u}^{\ -1}g]_0)\cdot ([g \overline{v^{-1}}\,]_-)^\iota \, (g^\iota)^{-1} \,
([\overline{u}^{\ -1}g]_+)^{\iota} \ .
\end{equation}
for any $g\in G^{u,v}$.
Other formulas for $\psi^{u,v}$ are:
$$\psi^{u,v}(g) =
u([\overline{u}^{\ -1}g]_0)[(\overline{v} g^{\iota})^{-1}]_+ \, \overline{v} \,
([\overline{u}^{\ -1} g]_+)^{\iota}\ .
$$
$$\psi^{u,v}(g)= u([\overline{u}^{\ -1}g]_0)\cdot ([g \overline{v^{-1}}\,]_-)^\iota \, \overline{u^{-1}}^{\,-1} \,
[\overline{u^{\ -1}} ((g)^\iota)^{-1}]_-\ .
$$

Also $\psi^{u,v}$ is symmetric: $(\psi^{u,v})^{-1}=\psi^{v,u}$. In particular, for $u=v$ the twist $\psi^{v,v}$ is an involution on $G^{v,v}$.

\end{lemma}

\begin{proof} Clearly, for any $h\in H$ and $x\in U\overline u U$ we have
$$[\overline{u}^{\ -1}hx]_0=[(\overline{u}^{\ -1}h\overline{u})\cdot \overline{u}^{\ -1}hx]_0=(\overline{u}^{\ -1}h\overline{u})\cdot [\overline{u}^{\ -1}hx]_0=\overline{u}^{\ -1}h\overline{u}=u^{-1}(h) \ .$$
Therefore, taking $g=hx$, where $h\in H$ and $x\in L^{u,v}$,  and taking into the account (\ref{eq:psi-u,v-x prime}) and (\ref{eq:psi-u,v-x}), we obtain the desirable formulas.
\end{proof}

Theorem \ref{th:t-through-x} admits the following obvious generalization.

\begin{theorem}
\label{th:t-through-x-h}
Let $\ii = (i_1, \ldots, i_m)$ be
a double reduced word for $(u,v)$, and suppose an element
$x \in G^{u,v}$ can be factored as
$x= h x_{i_1} (t_1) \cdots x_{i_m} (t_m)$,
with all $t_k$ nonzero elements of ${\mathcal F}$, and $h=diag(h_1,\ldots,h_n)\in H$.
Then the factorization parameters $h_1,\ldots, h_n$, $t_1,\ldots,t_m$
are determined by the following formulas:
\begin{equation}
\label{eq:h}
h_i =\Delta^{u^{-1}(i)}_{u,e}(x)
\end{equation}
for $i\in [1,n]$, and
\begin{equation}
\label{eq:t-through-x-h}
t_k =\begin{cases}
\Delta^i_{v_{<k} , u_{> k}}(y)^{-1}\Delta^i_{v_{<k} , u_{\geq k} } (y)=\Delta^{i+1}_{v_{<k}, u_{\geq k}}(y)^{-1}\Delta^{i+1}_{v_{<k} , u_{> k} } (y) & \text{if $i_k<0$}\\
\Delta^i_{v_{\le k}, u_{> k}} (y)^{-1}\Delta^{i+1}_{v_{<k} , u_{> k}} (y)
=\Delta^i_{v_{< k}, u_{> k}} (y)^{-1}\Delta^{i+1}_{v_{\le k} , u_{> k}} (y)  & \text{if $i_k>0$}
\end{cases}
\end{equation}
where $y=\psi^{u,v} (x)$ and $i=|i_k|$.
\end{theorem}

The following two special cases of Theorem \ref{th:t-through-x-h} will be of
particular importance: $(u,v) = (e, w_0)$
and $(u,v) = (w_), e)$ where $w_0$ is the longest element in $S_n$.
In these cases, Definition~\ref{def:twist} and Theorem~\ref{th:psi-regularity}
can be simplified as follows.

The formula  (\ref{eq:t-through-x})
now takes the following form.

\begin{corollary}
\label{cor:factors special}
Let $\ii = (i_1, \dots, i_m)$ be a reduced word for $w\in S_n$, and  $t_1, \dots, t_m$
be non-zero elements of ${\mathcal F}$.

\smallskip

\noindent (i)  If $x = x_{i_1} (t_1) \cdots x_{i_m} (t_m)$ then the factorization parameters $h_1,\ldots,h_n$ and
$t_1,\ldots,t_m$ are given by
$$h_i =\Delta^{i}_{e,e}(x)=x_{ii} $$
for $i\in [1,n]$, and
$$t_k =
\Delta^i_{s_{i_1} \cdots s_{i_k}, e} (y)^{-1}\Delta^{i+1}_{s_{i_1} \cdots s_{i_{k-1}} , e} (y)
=\Delta^i_{s_{i_1} \cdots s_{i_{k-1}}, e} (y)^{-1}\Delta^{i+1}_{s_{i_1} \cdots s_{i_k} , e} (y)
\ , $$
where $y=\psi^{e,w}(x)$ is given by {\rm (\ref{eq:psi e-wnot})}, and $i=i_k$.

\smallskip

\noindent (ii) If $x = hx_{- i_1} (t_1) \cdots x_{- i_m} (t_m)$ then the factorization parameters $h_1,\ldots,h_n$ and
$t_1,\ldots,t_m$ are given by
$$h_i =\Delta^{w^{-1}(i)}_{w,e}(x) $$
for $i\in [1,n]$, and
$$t_k =  \Delta^i_{e, s_{i_m} \cdots s_{i_{k+1}}} (y)^{-1} \Delta^i_{e, s_{i_m} \cdots s_{i_k} } (y)
=\Delta^{i+1}_{e, s_{i_m} \cdots s_{i_k}} (y)^{-1} \Delta^{i+1}_{e, s_{i_m} \cdots s_{i_{k+1}} } (y)\ , $$
where $y=\psi^{w,e}(x)$ is given by {\rm (\ref{eq:psi wnot-e})}.

\end{corollary}

\section{Other factorizations in $GL_n({\mathcal F})$ and the maximal twist $\psi^{\wnot ,\wnot }$}
\label{sec:other factorizations}



%
%
%
%
%
%

In this section we will provide some explicit factorizations in $G^{u,\wnot }$ and $G^{\wnot ,v}$. Let us consider a factorization of $x\in G^{u,\wnot }$ of the form:
\begin{equation}
\label{eq:positive standard}
x=x_-\cdot x^{(n-1)}x^{(n-2)}\cdots x^{(1)}\ ,
\end{equation}
where $x^-\in G^{u,e}$ and $x^{(m)}\in L^{e,s_m s_{m+1}\cdots s_{n-1}}$ is given by:
$$x^{(m)}=x_m(t_{m,m})x_{m+1}(t_{m,m+1})\cdots x_{n-1}(t_{m,n-1})$$
for $m\in [1,n-1]$.

\begin{lemma} In the notation of (\ref{eq:positive standard}), we have:
$$t_{m,k}=\Delta^{m,k}_{[1,m],[k-m+1,k]}(x)^{-1}\Delta^{m,k+1}_{[1,m],[k-m+2,k+1]}(x)$$
for all $1\le m\le k\le n-1$.

\end{lemma}

\begin{proof} Follows immediately from Theorem \ref{th:upper factorization}.
\end{proof}

\begin{lemma}
In the notation of (\ref{eq:positive standard}), we have:

$$t_{ij}=\Delta^{i,j}_{[1,i]\cup [n+i+1-j,n], [1,j]} (y)^{-1}\Delta^{i,j+1}_{[1,i]\cup [n+i-j,n] , [1,j+1]} (y)$$
for all $1\le i\le j< n$, where  $y=\psi^{u,\wnot }(x)$.
\end{lemma}

\begin{proof} Denote by $\ii_0$ the following {\it standard} reduced word for $\wnot $:
$$\ii_0=(n-1;n-2,n-1;\ldots;1,2,\ldots,n-1)\ .$$
It is convenient to identify $\ii_0$ with the sequence of pairs:
$$(n-1,n-1);(n-2,n-2),(n-2,n-1);\ldots;(1,1),(1,2),\ldots,(1,n-1)\ .$$

Let $\ii_-$ be any reduced word for $u\in S_n$. Then we put $\ii_-$ and $\ii_0$ into a separated word $\ii=(\ii_-,\ii_0)$ for the element $(u,\wnot )\in S_n\times S_n$.
Denote by $\wnot ^{(i,n)}$ the longest element of the subgroup of $S_n$ generated by the simple transpositions $s_i,s_{i+1},\ldots,s_{n-1}$.

Then in the notation of (\ref{eq:vlessk}) we have for the position $k$ of $\ii$ corresponding to the pair $(i,j)$:

\noindent $v_{\le k}=\wnot ^{(i+1,n)}s_is_{i+1}\cdots s_{j}, v_{<k}=\wnot ^{(i+1,n)}s_is_{i+1}\cdots s_{j-1}$,

\medskip


\noindent $v_{\le k}(j+1)=\wnot ^{(i+1,n)}s_is_{i+1}\cdots s_{j}(j+1)=\wnot ^{(i+1,n)}(i)=i$

\noindent $v_{<k}(j)=\wnot ^{(i+1,n)}s_is_{i+1}\cdots s_{j-1}(j)=\wnot ^{(i+1,n)}(i)=i$



\noindent $v_{\le k}[1,j+1]=\wnot ^{(i+1,n)}s_is_{i+1}\cdots s_{j-1}[1,j+1]=\wnot ^{(i+1,n)}[1,j+1]=[1,i]\cup [n+i-j,n]$

\noindent $v_{< k}[1,j]=\wnot ^{(i+1,n)}s_is_{i+1}\cdots s_{j-1}[1,j]=\wnot ^{(i+1,n)}[1,j]=[1,i]\cup [n+i+1-j,n]$


On the other  hand, taking  (\ref{eq:t-through-x-h}) for $\ii=(\ii_-,\ii_0)$ with $i_k=j$, yields the following formula
$$t_k=\Delta^j_{v_{< k}, e} (y)^{-1}\Delta^{j+1}_{v_{\le k} , e} (y)$$
which, after substituting the results of the above computations, implies the desirable formula for $t_k=t_{ij}$.

The lemma is proved.
\end{proof}

The above facts imply an immediate corollary.

\begin{corollary} For any $u\in S_n$ the twist map $\psi^{u,\wnot }$ satisfies:
$$\Delta^{i,j}_{[1,i]\cup [n+i+1-j,n], [1,j]} (\psi^{u,\wnot }(x))^{-1}\Delta^{i,j+1}_{[1,i]\cup [n+i-j,n] , [1,j+1]} (\psi^{u,\wnot }(x))=$$
$$=\Delta^{i,j}_{[1,i],[j+1-i,j]}(x)^{-1}\Delta^{i,j+1}_{[1,i],[j-i+2,j+1]}(x)$$
for all $1\le i\le j\le n-1$.

\end{corollary}

Let us consider a factorization of $x\in G^{\wnot ,v}$ of the form
\begin{equation}
\label{eq:negative standard}
x_-=h \cdot x_-^{(n-1)}x_-^{(n-2)}\cdots x_-^{(1)}\cdot x_+
\end{equation}
where $x_+\in L^{e,v}$, $h\in H$,  and $x_-^{(m)}\in L^{s_m\cdots s_{n-1}s_m,e}$ is of the form:
$$x_-^{(m)}=x_{-m}(\tau_{m,m})x_{-(m-1)}(\tau_{m,m+1})\cdots x_{-(n-1)}(\tau_{m,n-1})$$
for $m\in [1,n-1]$.

The following result generalizes  the factorization from Section
3.2.

\begin{proposition} In the notation of (\ref{eq:negative standard}) we have
$$h_n=x_{n1}, h_{n-1}=-\left |\begin{matrix}
x_{n-1,1}&\boxed {x_{n-1,2}}\\
x_{n,1}&x_{n, 2}\end{matrix} \right |,
\dots , h_1=(-1)^{n-1}|x|_{1n}\ ,$$
that is,
$$h_m=\Delta^{m,n+1-m}_{[m,n],[1,n+1-m]}(x)$$
for $m\in [1,n]$,

and
$$\tau_{m,k}=(-1)^{k-m}\left |\begin{matrix}
x_{m,1}&\dots &\boxed {x_{m, k+1-m}}\\
            &\dots &                       \\
x_{k,1}&\dots &x_{k,k+1-m}\end{matrix} \right |^{-1}
h_m$$
for all $1\le m\le k <n$,
i.e.,
$$\tau_{m,k}=\Delta^{m,k+1-m}_{[m,k],[1,k+1-m]}(x)^{-1}\Delta^{m,n+1-m}_{[m,n],[1,n+1-m]}(x) \ .$$


\end{proposition}

\bigskip
The proof is similar to the proof of Theorem \ref{th:upper factorization}.

\begin{example} Let $n=3$. Then in the factorization
$$x=h\cdot x_{-2}(\tau_{22})x_{-1}(\tau_{11})x_{-2}(\tau_{12})\cdot x_+ \ ,$$
where $h\in H$ and $x_+\in U$, we have:
$$\tau_{11}=x_{11}^{-1}\Delta^{1,3}_{123,123}(x),~\tau_{12}=\Delta^{1,2}_{12,12}(x)^{-1}\Delta^{1,3}_{123,123}(x),~\tau_{22}=x_{21}^{-1}\Delta^{2,1}_{23,12}(x) \ .$$

\end{example}
Our next result is a direct consequence of Theorem \ref{eq:t-through-x-h}.

\begin{lemma}
In the notation of (\ref{eq:negative standard}), we have:
$$\tau_{ij}=\Delta^{j,j+1-i}_{[1,j], [n+2-i,n]\cup [1,j+1-i]}(y)^{-1}\Delta^{j,n+1-i}_{[1,j] , [n+1-i,n]\cup [1,j-i] } (y)$$
for all $1\le i\le j< n$, where  $y=\psi^{\wnot ,e}(x)$.
\end{lemma}

\begin{proof} Recall that $\ii_0=(n-1;n-2,n-1;\ldots;1,2,\ldots,n-1)$
is the standard reduced word for $\wnot $ and that we conveniently identified
$\ii_0$ with the sequence of pairs:
$$(n-1,n-1);(n-2,n-2),(n-2,n-1);\ldots;(1,1),(1,2),\ldots,(1,n-1)\ .$$

Let $\ii_+$ be any reduced word for $v\in S_n$. Then we put $-\ii_0$ and $\ii_+$ into a separated word $\ii=(-\ii_0,\ii_+)$ for the element $(\wnot ,v)\in S_n\times S_n$.
Recall that $\wnot ^{(i,n)}$ denotes the longest element of the subgroup of $S_n$ generated by the simple transpositions $s_i,s_{i+1},\ldots,s_{n-1}$.

Then in the notation of (\ref{eq:ulessk}) we have for the position $k$ of $\ii$ corresponding to the pair $(i,j)$:

\noindent $u_{\ge k}=\wnot \wnot ^{(i,n)}s_{n-1}s_{n-2}\cdots s_{j}, u_{>k}=\wnot \wnot ^{(i,n)}s_{n-1}s_{n-2}\cdots s_{j+1}$,

\medskip

%
%
%

\noindent $u_{\ge k}(j)=\wnot \wnot ^{(i,n)}(n)=\wnot (i)=n+1-i$

\noindent $u_{>k}(j)=\wnot \wnot ^{(i,n)}(j)=\wnot (n+i-j)=j+1-i$

\noindent $u_{\ge  k}[1,j]=\wnot \wnot ^{(i,n)}([1,j-1]\cup \{n\})=\wnot ([1,i]\cup [n+i-j,n])=[n+1-i,n]\cup [1,j-i]$

\noindent $u_{> k}[1,j]=\wnot \wnot ^{(i,n)}([1,j])=\wnot ([1,i-1]\cup [n+i-j,n])=[n+2-i,n]\cup [1,j+1-i]$

On the other  hand, taking  (\ref{eq:t-through-x-h}) for $\ii=(-\ii_0,\ii_+)$ with $i_k=-j$, yields the following formula
$$t_k=\Delta^{j}_{e, u_{> k}}(y)^{-1}\Delta^{j}_{e , u_{\ge k} } (y) \ ,$$
which, after substituting the results of the above computations, implies the desirable formula for $\tau_k=\tau_{ij}$.

The lemma is proved. \end{proof}

The above facts imply an immediate corollary.

\begin{corollary} For any $v\in S_n$ the twist map $\psi^{\wnot ,v}$ satisfies:
$$\Delta^{j,j+1-i}_{[1,j], [n+2-i,n]\cup [1,j+1-i]}(\psi^{\wnot ,v }(x))^{-1}\Delta^{j,n+1-i}_{[1,j] , [n+1-i,n]\cup [1,j-i] } (\psi^{\wnot ,v }(x))=$$
$$\Delta^{i,j+1-i}_{[i,j],[1,j+1-i]}(x)^{-1}\Delta^{i,n+1-i}_{[i,n],[1,n+1-i]}(x)$$
for all $1\le i\le j\le n-1$.

\end{corollary}

The above results allow us to completely compute the twist $\psi^{\wnot ,\wnot }$ in terms of positive quasiminors.

\begin{theorem}
\label{th:double ratios}
For each $x\in G$ we have (with the notation $y=\psi^{\wnot ,\wnot }(x)$):

$$\Delta^{n+1-i,i}_{[n+1-i,n],[1,i]}(y)=\Delta^{n+1-i,i}_{[n+1-i,n],[1,i]}(x)$$
for $i\in [1,n]$, and:
$$\Delta^{i,j}_{[1,i]\cup [n+i+1-j,n], [1,j]} (y)^{-1}\Delta^{i,j+1}_{[1,i]
\cup [n+i-j,n] , [1,j+1]} (y)=
\Delta^{i,j}_{[1,i],[j+1-i,j]}(x)^{-1}\Delta^{i,j+1}_{[1,i],[j-i+2,j+1]}(x),$$
\smallskip
$$\Delta^{i,j}_{[1,i],[j+1-i,j]}(y)^{-1}\Delta^{i,j+1}_{[1,i],[j-i+2,j+1]}(y)=
\Delta^{i,j}_{[1,i]\cup [n+i+1-j,n], [1,j]} (x)^{-1}\Delta^{i,j+1}_{[1,i]
\cup [n+i-j,n] , [1,j+1]} (x),$$
\smallskip
$$\Delta^{j,j+1-i}_{[1,j], [n+2-i,n]\cup [1,j+1-i]}(y)^{-1}
\Delta^{j,n+1-i}_{[1,j] , [n+1-i,n]\cup [1,j-i]} (y)=
\Delta^{i,j+1-i}_{[i,j],[1,j+1-i]}(x)^{-1}
\Delta^{i,n+1-i}_{[i,n],[1,n+1-i]}(x),$$
\smallskip
$$\Delta^{i,j+1-i}_{[i,j],[1,j+1-i]}(y)^{-1}
\Delta^{i,n+1-i}_{[i,n],[1,n+1-i]}(y)=
\Delta^{j,j+1-i}_{[1,j], [n+2-i,n]\cup [1,j+1-i]}(x)^{-1}
\Delta^{j,n+1-i}_{[1,j] , [n+1-i,n]\cup [1,j-i] } (x)$$
for all $1\le i\le j\le n-1$.

\end{theorem}

The above result allows to compute explicitly a large number of
positive quasiminors for maximally twisted matrices and to get
other relations.

\begin{corollary}

In the notation of Theorem \ref{th:double ratios} we have

$$
\Delta ^{i, j+1-i}_{[i,j], [1,j+1-i]}(y)=
\Delta ^{i, n+1-i}_{[i,n], [1,n+1-i]}(x)
\Delta ^{j, n+1-i}_{[1,j], [n+1-i, n]\cup [1,j-i]}(x)^{-1}
\Delta ^{j, j+1-i}_{[1,j], [n+2-i, n]\cup [1,j+1-i}(x)
$$
for all $1\le i\le j\le n-1$.

Also,
$$
\Delta ^{i,i}_{[1,i], [1,i]}(y)^{-1} \Delta ^{i, j}_{[1,i]\cup
[n+i-j+1, n], [1,j]}(y)= \Delta ^{i,i}_{[1,i], [1,i]}(x)^{-1}
\Delta ^{i, j}_{[1,i], [j-i+1,j]}(y),
$$
$$
\Delta ^{i,i}_{[1,i], [1,i]}(y)^{-1} \Delta ^{i, j}_{[1,i],
[j-i+1,j]}(y)= \Delta ^{i,i}_{[1,i], [1,i]}(x)^{-1} \Delta ^{i,
j}_{[1,i]\cup [n+i-j+1], [j-i+1,j]}(y).
$$
\end{corollary}

\end{document}